\documentclass[12pt,reqno]{amsart}
\usepackage{exscale}
\usepackage{relsize}
\usepackage{amsfonts}
\usepackage{a4wide}
\usepackage{mathrsfs,amsmath,makecell,amssymb}
\usepackage{makecell}
\usepackage{amssymb}
\usepackage{multicol}
\usepackage{enumitem}
\usepackage{color}
\usepackage[colorlinks=true,urlcolor=blue,citecolor=red,linkcolor=blue,linktocpage,pdfpagelabels,bookmarksnumbered,bookmarksopen]{hyperref}
\hypersetup{urlcolor=blue, citecolor=red, linkcolor=blue}

\allowdisplaybreaks
\numberwithin{equation}{section}
\newtheorem{theorem}{Theorem}[section]
\newtheorem{proposition}[theorem]{Proposition}

\newtheorem{lemma}[theorem]{Lemma}

\theoremstyle{definition}

\newtheorem{remark}[theorem]{Remark}

\newcommand{\R}{\mathbb{R}}
\newcommand{\dis}{\displaystyle}

\newcommand{\eps}{\varepsilon}

\def\R{\mathbb{R}}

\def\E{\mathcal{E}}

\def\Qc{\mathcal{Q}}

\def\X{\mathcal{X}}
\def\C{\mathcal{C}}

\def\RN{\mathbb{R}^N}

\def\n{\nabla}
\def\de{\partial}
\def\a{\alpha}

\def\l{\lambda}
\def\g{\gamma}
\def\s{\sigma}
\def\t{\theta}
\def\irn{\int_{\RN}}
\def\dis{\displaystyle}
\newcommand{\weakto}{\rightharpoonup}

\newcommand{\cC}{{\mathcal C}}

\usepackage[colorinlistoftodos]{todonotes}
\makeatletter
\providecommand\@dotsep{5}
\def\listtodoname{List of Todos}
\def\listoftodos{\@starttoc{tdo}\listtodoname}
\makeatother
\allowdisplaybreaks

\begin{document}

\title[Normalized solutions to Born-Infeld and quasilinear problems]
{Normalized solutions to Born-Infeld\\ and quasilinear problems}

\author[L. Baldelli]{Laura Baldelli} \email{labaldelli@ugr.es}

\author[J. Mederski]{Jarosław Mederski} \email{jmederski@impan.pl}

\author[A. Pomponio]{Alessio Pomponio}
\email{alessio.pomponio@poliba.it}

\address[Baldelli]{IMAG, Departamento de Análisis Matemático -- Universidad de Granada -- Campus Fuentenueva -- 18071 Granada, Spain}

\address[Mederski]{Institute of Mathematics -- Polish Academy of Sciences -- ul. \'Sniadeckich 8 -- 00-656 Warsaw, Poland}

\address[Pomponio]{Dipartimento di Meccanica, Matematica e Management -- Politecnico di Bari -- Via Orabona 4 -- 70125 Bari, Italy}

\begin{abstract}
The paper concerns the existence of normalized solutions to a large class of quasilinear problems, including the well-known Born-Infeld operator. In the mass subcritical cases, we study a global minimization problem and obtain a ground state solution for a $(2,q)$-type operator which implies the existence of solutions to the Born-Infeld problem.
We also deal with the mass critical and mass supercritical cases for quasilinear problems involving the $(2,q)$-type operator.
\end{abstract}


\maketitle

\section{Introduction}

The electromagnetic theory announced by Born and Infeld, see \cite{Bnat,B,BInat,BI}, is a nonlinear alternative to the classical Maxwell theory. The importance of this theory lies in giving a unitarian point of view to describe electrodynamics and in notable feature to be a fine answer to the well-known {\em infinite-energy problem}: indeed, the electromagnetic field generated by a point charge has finite energy in the Born-Infeld model in which a central role is played by the following peculiar differential operator
\begin{equation}\label{biop}
\Qc(u)=-{\rm div}\left(q(|\n u|^2)\n u \right)\qquad\text{where} \qquad q(s)=\frac{1}{\sqrt{1-s}}.
\end{equation}
The mean curvature operator in Lorentz-Minkowski space can be described by this kind of operator, see for instance \cite{BS,CY}. 

The problem 
\[
-{\rm div}\left(\dfrac{\n u}{\sqrt{1-|\n u|^2}}\right)=\varrho, \qquad\hbox{ in }\RN,
\]
with $\varrho$ assigned, was studied by many authors, see \cite{APS, BCP, BCF,BDP,BDPR,Bon-Iac2,Bon-Iac3,BIMM,H,K,K-corr}, under different assumptions about $\varrho$. 

Few is still known, at contrary, when $\varrho$ is replaced by a nonlinearity, namely for an equation of the type
\begin{equation}\label{cong}
-{\rm div}\left(\dfrac{\n u}{\sqrt{1-|\n u|^2}}\right)=g(u), \qquad\hbox{ in }\RN.
\end{equation}
Here classical variational techniques do not work directly for this problem, due to the particular nature of the operator $\Qc$. The related functional
\[
J(u)=\irn \left(1-\sqrt{1-|\n u|^2}\right) -\irn G(u)\, dx,
\]
where $G$ is a primitive of $g$, is well defined only 
under the condition $|\n u|\le 1$ a.e. in $\RN$ and it is not regular on the set $\{x\in \RN :|\n u|=1\}$: this requires different and non-standard strategies. 
As far as we know, a first attempt in this direction  using variational methods was made in  \cite{BDD}, where $g(s)=|s|^{p-2}s$, for $p>2^*=\frac{2N}{N-2}$ and $N\ge 3$. By means of suitable truncation arguments (that will be crucial in our approach, as we will see later), the existence of {\em finite energy} solutions is proved. Later on in \cite{MP}, this strategy has been extended to treat a larger class of nonlinearities. 
\\
We mention, moreover, \cite{A,A2,P} where \eqref{cong} has been studied by means of ODE-techniques finding solutions which could have infinite energy.\\
\indent Our principal aim in this paper is to study {\em normalized solutions} to \eqref{cong} and related quasilinear problems, that is solutions such that the quantity, the so-called {\em mass} 
\begin{equation}\label{eq:mass}
\int_{\R^N} |u|^2 \, dx =\rho^2
\end{equation}
is prescribed for given $\rho>0$. In the various physical contexts this quantity represents e.g., the power supply in electromagnetism and the total number of atoms in the Bose--Einstein condensation. Moreover the mass arises when one looks for standing-wave solutions to the evolution equations and the mass is one of the conservation laws. When the left hand side of \eqref{cong} is replaced by a homogenous operator and $g(s)=|s|^{p-2}s-s$, then any nontrivial solution can be rescaled to the normalized one. Indeed, recall that Kwong \cite{Kwong} showed that there exists a unique positive radially symmetric solution to the equation
$$-\Delta u + u=|u|^{p-2}u$$
and up to some rescaling $\alpha u(\gamma \cdot)$ for some $\alpha,\gamma>0$ is the unique positive radially symmetric solution  $-\Delta u +\lambda u=|u|^{p-2}u$ for some $\lambda>0$ with the prescribed mass \eqref{eq:mass}.  The classical results in this area go back to the works \cite{CazenaveLions,s80, s82}, and the most recent ones are presented in the works e.g. \cite{BartschSoaveJFA,BartschSoaveJFACorr,BiegMed,J7,JeanjeanLu}, see also the references therein.

We would like to emphasize that the rescaling argument is not available for the nonlinear operator $\Qc$, moreover in order to find a solution to \eqref{cong} we have to keep the $L^\infty$-upper bound on $\nabla u$, which is an additional difficulty. Up to our knowledge, there is no result in the literature on the existence of normalized solution in presence of Born-Infeld operator. 
In the present paper we are interested in normalized solutions to even a larger class of operators, including the Born-Infeld operator, with power-type nonlinearity, of the form
\begin{equation}\label{P}\tag{${\mathcal P}_a$}
\begin{cases}
-\mbox{div}(a(|\nabla u|^2)\nabla u)=|u|^{p-2}u-\lambda u, \qquad & \text{in}\,\,\mathbb R^N,\\
u(x)\to 0, & |x|\to+\infty,\\
\int_{\R^N}|u|^2\,dx=\rho^2,
\end{cases}\end{equation}
where  $N\ge 3$, $2<p<2^*$ and under the following assumptions on $a$
\begin{enumerate}[label=(a\arabic{*}),ref=a\arabic{*}]
\setcounter{enumi}{-1}
	\item \label{a0}$a:[0,1) \to (0,+\infty)$ is continuous, of class $\mathcal{C}^1$ on $(0,1)$;
	\item \label{a1} 
	$\displaystyle\lim_{s\to 1^-}a(s)=+\infty;$
	\item \label{a2} $a$ is  increasing.
\end{enumerate}

Considering $a(s)=(1-s)^\alpha$ with $\alpha<0$, it satisfies \eqref{a0}--\eqref{a2}. Moreover, we get the operator $\Qc$ for $\alpha=-1/2$. Another important example is the following general mean curvature operator arising in the study of hypersurfaces in the Lorentz–Minkowski space $\mathbb L^{N+1}$ and in $\mathbb R^{N+1}$
$$a(s)=\beta (1-s)^{-1/2}-\gamma (1+s)^{-1/2},\quad \beta>0, \gamma\geq 0,$$
see \cite{CY, MP} and references therein.

The first aim of the paper is to find normalized weak solutions in 
\begin{equation}\label{Sc}
S_{\rho}: =\Big\{ u \in \X: \int_{{\mathbb{R}^N}} {|u|}^2dx=\rho^2  \Big\},
\end{equation}
where $\X=H^1(\R^N)\cap D^{1,q}(\R^N)$ for some  $q>1$ defined below, and we will study the {\em $L^2$-subcritical case}, i.e. $p<2+4/N$, and the {\em $L^2$-supercritical case}, i.e. $p>2+4/N$, by using different techniques. 

Concerning the Born-Infeld problem \eqref{P} we get the following result

\begin{theorem}\label{main}
Let $\rho>0$ and $2<p<2+\frac{4}{N}$. Then, there exists a nonnegative radial solution $u\in\X$ to \eqref{P}. 
\end{theorem}

It remains as an open question whether there are solutions to \eqref{P} in the $L^2$-supercritical case.

The outline of the proof will be the following. Since $a$ is diverging near $1$ and so the differential operator is singular, following \cite{BDD,MP}, we truncate the operator and, as an intermediate step, we study the existence of normalized solutions for the truncated problem. 

In the proof of Theorem \ref{main}, the study of solutions to the following problem, which is of some interest by itself, will be crucial.  Namely, for
 $N\ge 3$, we consider the following {\em quasilinear problem}
\begin{equation}\label{Pb}\tag{${\mathcal P}_b$}
\begin{cases}
-\mbox{div}(b(|\nabla u|^2)\nabla u)=|u|^{p-2}u-\lambda u, \qquad & \text{in}\,\,\mathbb R^N,\\
u(x)\to 0, & |x|\to+\infty,\\
\int_{\R^N}|u|^2\,dx=\rho^2,
\end{cases}\end{equation}
under the assumptions on $b$:
\begin{enumerate}[label=(b\arabic{*}),ref=b\arabic{*}]
\setcounter{enumi}{-1}
	\item \label{b0}$b:[0,+\infty) \to (0,+\infty)$ is of class $\cC^1$;
	\item \label{b1} there exist $q\in \left(1,+\infty\right)$ and $c_1,c_2>0$ such that
\begin{align}\label{growthb1}
c_1(s^2+|s|^q)&\le b(s^2)s^2\le c_2(s^2+|s|^q)
\\
\label{growthb2}
c_1(s^2+|s|^q)&\le B(s^2)\le c_2(s^2+|s|^q)
\end{align}
where $B(s)=\int_0^s b(t) dt$;
	\item \label{b2} $b$ is increasing.
\end{enumerate}
Firstly we consider the  $L^2$-subcritical case and $L^q$-subcritical case, that is 
\begin{equation}\label{h1.1}
2<p< \biggl(1+\frac{2}{N}\biggl)\min\{2,q\}.
\end{equation}

We obtain a {\em ground state solution}, that is a minimal energy among all solutions that belong to $S_\rho$. More precisely we show that there exists a global minimizer for
$$ m(\rho)\!:=\!\inf _{u \in S_\rho} I(u),$$
where 
\begin{equation*}
I(u)= \frac{1}{2}\int_{\mathbb R^N} B(|\nabla u|^2) \, dx-\frac{1}{p}{\|u \|}_p^p.
\end{equation*}

\begin{theorem}\label{mainb}
Assume $\rho>0$, \eqref{b0}--\eqref{b2} and \eqref{h1.1}. Then, there exists a nonnegative radially decreasing normalized ground state solution $u_\rho\in\X$ to \eqref{Pb}. In particular, $I(u_\rho)=m(\rho)<0$.
\end{theorem}

Taking inspiration from classical Cazenave, Lions \cite{CazenaveLions}, and Stuart \cite{s82}, we prove the existence of a radial normalized solution to \eqref{Pb} by showing that the global minimization problem is attained by excluding vanishing and dichotomy and overcoming the nonhomogenity of the operator taken into account, see also recent works \cite{BY,JeanjeanLu} and references therein. Recall that Theorem \ref{mainb} contains the result of \cite{BY}, where the $(2,q)$-Laplacian operator, i.e. $b(s)=1+s^{q-2}$ has been considered. A more general class $b$ plays a crucial role due to the application to the Born-Infeld problem \eqref{P} and makes the variational approach more difficult as we shall see later.

Following \cite{BDD,MP}, by means of Theorem \ref{mainb}, we are able to obtain a radial normalized solution for a truncated version of problem \eqref{P} and,  by proving uniformly boundedness of the gradient of such solution, we manage to turn it into a solution of the main problem \eqref{P} which includes also the Born-Infeld equation.
One of the novelties of the paper consists of analysing the constrained Born-Infeld problem, which required different tools with respect to the unconstrained case. For instance, in treating Born-Infeld type equations we need to treat a suitable truncated problem, slightly different with respect to the unconstrained case characterized by a $\cC^1$-function. Lastly, we assume \eqref{a2}, which is a weaker assumption than  the strict convexity assumed in \cite{MP}.

If \eqref{h1.1} is violated, we can find solutions in the following $L^2$-supercritical case and $L^q$-supercritical case case, that is we assume
\begin{equation}\label{eq:L2super}
\biggl(1+\frac{2}{N}\biggl)\max\{2,q\}<p<2^*,\; q>\frac{2N}{N+2}
\end{equation}
and we consider also the following assumption on $b$
\begin{enumerate}[label=(b\arabic{*}),ref=b\arabic{*}]
	\setcounter{enumi}{2}
	\item \label{b3}$\a b(s)s\le B(s)<\frac{2(N+p)}{Np}b(s)s  \text{ for any } s>0\hbox{ and some }\a>\frac{2(N+2)}{Np}$.
\end{enumerate}

\begin{theorem}\label{mainb2}
	Assume $\rho>0$, \eqref{b0}--\eqref{b3} and \eqref{eq:L2super}. Then, there exists a radial normalized solution $u\in\X$ to \eqref{Pb}.
\end{theorem}

We will employ a mountain pass approach to
prove Theorem \ref{mainb2}. In particular, we will observe that the related functional possesses a mountain pass geometry and it is possible to define a mountain pass
minimax value. Then, by
Ekeland’s principle, we can find a Palais–Smale sequence but the lack of the Palais–Smale compactness condition is difficult to overcome. A key of our argument is to find a Palais–Smale sequence with an extra property related to Pohozaev's identity: for this purpose we introduce an auxiliary functional taking inspiration by \cite{HIT,J7}.

Finally, we are able to give some partial nonexistence results for normalized ground state solutions for \eqref{Pb}, in the $L^2$-  and $L^q$-mass critical cases. More precisely, the following holds.

\begin{theorem}\label{th1.2}
Assume $N\ge 3$, $p=2+4/N$, \eqref{b0}--\eqref{b1}. Then there exist $0<\rho_*<\rho^*$ such that  
\begin{enumerate}[label=(\roman{*}),ref=\roman{*}]
\item \label{i-th1.2} if $q>1$ and $0\!<\!\rho\!\leq\! \rho_*$, then $m(\rho)\!=\!0$ and there is no minimizer of $m(\rho)$; 
\item  \label{ii-th1.2}if $1\!<\!q\!<\!2$ and $\rho\!>\!\rho^*$, then $m(\rho)\!=\!-\infty$ and there is no minimizer of $m(\rho)$.
\end{enumerate}
\end{theorem}

\begin{theorem}\label{th1.2bis}
Assume $N\ge 3$, $p\!=q+2q/N$, \eqref{b0}--\eqref{b1}. Then there exist $0<\hat \rho_{*}<\hat{\rho}^{*}$ such that   
\begin{enumerate}[label=(\roman{*}),ref=\roman{*}]
\item \label{i-th1.2bis} if $1<q<N$ and $0\!<\!\rho\!\leq\! \hat \rho_{*}$, then $m(\rho)\!=\!0$ and there is no minimizer of $m(\rho)$;  
\item  \label{ii-th1.2bis} if $2\!<\!q\!<\!N$ and $\rho\!>\!\hat{\rho}^{*}$, then $m(\rho)\!=\!-\infty$ and there is no minimizer of $m(\rho)$.
\end{enumerate}
\end{theorem}

The paper is organized as follows. In Section \ref{prel}, we introduce our functional framework and some technical results. Sections \ref{app} and \ref{res} are focused on proving the existence of solutions in the subcritical case to \eqref{Pb} and \eqref{P}, respectively. 
The supercritical case for the generalized operator is treated in Section \ref{super}, together with the difficulties in considering the Born-Infeld problem under this light. Finally, Section \ref{crit} is devoted to the critical cases.

\section{Preliminaries}\label{prel}

In this section we state some preliminary features, starting from the functional setting. 
 In the following, we fix  $q>1$ and $N\ge 3$. 
Let
$\X$ be the completion of $\cC_0^{\infty}(\R^N)$ with respect to the norm
$$\|u\|_\X=\big(\|\nabla u\|_2^2+\|u\|_2^2+\|\nabla u\|_q^2\big)^{1/2}.$$ 
Here and in what follows $\|\cdot\|_p$ stands for $L^p$-norm. Moreover $B_R(0)$ denotes and open ball centred at $0$ with radius $R>0$, $ o(1)$ is standard asymptotic notation.
Clearly  $\X=H^1(\R^N)\cap D^{1,q}(\R^N)\subset D^{1,2}(\R^N)\cap D^{1,q}(\R^N)$ equipped with the norm 
$$\|u\|_{2,q}=(\|\nabla u\|_2^2+\|\nabla u\|_q^2)^{1/2},$$
which is equivalent to the usual norm $\|\nabla u\|_2+\|\nabla u\|_q$. Moreover, since 
\[
\X \hookrightarrow   H^1(\R^N)\hookrightarrow L^s(\R^N),\qquad\text{for }
s\in [2,2^*]
\]
\[
\X \hookrightarrow D^{1,2}(\R^N)\cap D^{1,q}(\R^N)\hookrightarrow L^{s}(\R^N),\qquad\text{for }
s\in \begin{cases}
[q^*,2^*]& \text{if}\,\, 1<q\le 2\\
[2^*,q^*]& \text{if}\,\, 2<q<N\\
[2^*,+\infty)& \text{if}\,\, 2<q=N\\
[2^*,+\infty]& \text{if}\,\, q>N
\end{cases}
\]
we have that 
\[
\X \hookrightarrow    L^s(\R^N),\qquad\text{for }
s\in \begin{cases}
[2,2^*]& \text{if}\,\, 1<q\le2<N\\
[2,q^*]& \text{if}\,\, 2<q<N\\
[2,+\infty)& \text{if}\,\,  2<q=N\\
[2,+\infty]& \text{if}\,\, q>N
\end{cases}.
\]
In Section \ref{res}, we will need to restrict our consideration to the radial setting. So we denote 
\begin{align*}
\X_r:=\big\{u\in \X: u \hbox{ radially symmetric}\big\},\qquad
S_{\rho,r}:=\left\{u\in \X_r: \irn |u|^2 dx=\rho^2\right\},
\end{align*}
for $c>0$, and we have
 \[
\X_r  \hookrightarrow H^1_r(\R^N) \hookrightarrow \hookrightarrow L^s(\R^N),\qquad\text{for }
s\in \begin{cases}
(2,2^*)& \text{if}\,\, 1<q\le2\\
(2,q^*)& \text{if}\,\, 2<q<N\\
(2,+\infty)& \text{if}\,\, q\ge N
\end{cases}
\]
 see e.g. \cite{BDD,AW}.



The following lemma is a variant  of the well-known Strauss Lemma (see \cite[Lemma 2.3]{MP}, cf. \cite{S,SWW}).


\begin{lemma}\label{strauss}
	Let $N\ge 2$ and $q>N$. There exists $C=C(N,q) > 0$  such that for all $u \in \X$, there holds
	\[
	|u(x)| \le C|x|^{-\frac{N-1}2} \|u\|_{\X}, 
	\]
	for all $|x|\ge 1$.
\end{lemma}

Now, we report for completeness the Gagliardo-Nirenberg inequality, which will be very useful in what follows.

\begin{lemma} (Gagliardo-Nirenberg inequality, \cite{Wein}) \label{lem2.1}
Let $p\!\in\!(2,2^*)$ and $\delta_p\!=\!\frac{N(p-2)}{2p}\in(0,1)$. Then there exists a constant $\mathcal{C}_{N,p}\!=\!\left( \frac{p}{2 \|W_p\|^{p-2}_{2}} \right)^{\frac{1}{p}}\!>\!0$ such that
\begin{equation} \label{equ2.2}
 \|u\|_{p} \leq \mathcal{C}_{N,p} \left\|\nabla u\right\|_{2}^{\delta_p} \left\|u\right\|_{2}^{(1-\delta_p)}, \qquad \forall u \in {H}^{1}(\mathbb{R}^{N}),
\end{equation}
where $W_p$ is the unique positive radial solution of 
\begin{equation}\label{eqlem2.1}
-\Delta W\!+\!(\frac{1}{\delta_p}-\!1)W \!=\!\frac{2}{p\delta_p}|W|^{p-2}W.
\end{equation}
\end{lemma}

Moreover, to deal with the $L^q$-critical case we need to recall also the $L^q$-Gagliardo-Nirenberg inequality

\begin{lemma} ($L^q$-Gagliardo-Nirenberg inequality, \cite[Theorem 2.1]{MaEh}) \label{lem2.2}
Let $N\ge 3$, $q\!\in\!\left(\frac{2N}{N+2},N\right)$, $p\!\in\!(2,q^*)$ and $\nu_{p,q}\!=\!\frac{Nq(p-2)}{p[Nq-2(N-q)]}$. Then there exists a constant $\mathcal{K}_{N,p}\!>\!0$ such that
\begin{equation*} \label{pq-equ2.2}
 \|u\|_{p} \leq \mathcal{K}_{N,p} \left\|\nabla u\right\|_{q}^{\nu_{p,q}} \left\|u\right\|_{2}^{(1-\nu_{p,q})}, \qquad \forall u \in {D}^{1,q}(\mathbb{R}^{N})\cap{L}^{2}(\mathbb{R}^{N}),
\end{equation*}
where 
$$\mathcal{K}_{N,p}=\biggl[\frac{K}{\frac{1}{q}\|\nabla W_{p,q}\|_q^q+\frac{1}{2}\|W_{p,q}\|_2^2}\biggr],$$
$$K=(Nq+pq-2N)\cdot\biggr[\frac{[2(Nq-p(N-q))]^{p(N-q)-Nq}}{[qN(p-2)]^{N(p-2)}}\biggr]^{1/[Nq+pq-2N]},$$
and $W_{p,q}$ is the unique nonnegative radial solution of the following equation
$$-\Delta_q W+W=\zeta|W|^{p-2}W,$$
where $\zeta=\|\nabla W\|_q^q+\|W\|_2^2$ is the Lagrangian multiplier.
\end{lemma}

We end this section with the following useful result, of which we report the proof based on \cite{Fig2011}.
\begin{lemma}\label{lemdis}
Assume $b$ satisfying \eqref{b0}--\eqref{b2}. Then there exists a constant $C>0$ such that the following inequality holds
$$\langle b(|x|^2)x-b(|y|^2)y, x-y \rangle\ge C\bigl[|x-y|^2+|x-y|^q\bigr]$$
for all $x,y\in\mathbb R^N$.
\end{lemma}
\begin{proof}
First, note that
$$\langle b(|x|^2)x-b(|y|^2)y, x-y \rangle=\sum_{j=1}^N \left(b(|x|^2)x_j-b(|y|^2)y_j\right)(x_j-y_j).$$
On the other hand, since $b\in C^1(\mathbb R^+)$ for every $\xi,z\in\mathbb R^N$ we have
\begin{equation}\label{2.12}
\begin{aligned}
\sum_{i,j=1}^N \frac{\partial}{\partial z_i} \left(b(|z|^2)z_j\right)\xi_i\xi_j&=\sum_{i,j=1}^N b(|z|^2) \delta_{ij} \xi_i\xi_j+2\sum_{i,j=1}^N b'(|z|^2)z_iz_j\xi_i\xi_j\\
&=b(|z|^2) |\xi|^2+2b'(|z|^2)\sum_{i,j=1}^N z_iz_j\xi_i\xi_j
\end{aligned}\end{equation}
Assume, without loss of generality, that $|y|\ge|x|$, thus $\frac{1}{2}|x-y|\le|y|$. Moreover, for $t\in[0,\frac{1}{4}]$ we get
\begin{equation}\label{xyt}
|y+t(x-y)|\ge |y|-t|x-y|\ge\frac{1}{4}|x-y|
\end{equation}
Taking $\xi=x-y$ and $z=y+t(x-y)$ and by using the Fundamental Theorem of Calculus in $[0,1]$ with $F_j(t)=b(|y+t(x-y)|^2)(y_j+t(x_j-y_j))$, we have
$$\sum_{j=1}^N \big(b(|x|^2)x_j-b(|y|^2)y_j\big)(x_j-y_j)=\int_0^1 \sum_{i,j=1}^N \frac{\partial}{\partial z_i} (b(|z|^2)z_j)\xi_i\xi_j dt$$
Now, from the equality above, using \eqref{2.12} it holds
\begin{equation}\label{dif}
\begin{aligned}
\langle b(|x|^2)x-b(|y|^2)y, x-y \rangle&= \int_0^1 \biggl(b(|z|^2) |\xi|^2+2b'(|z|^2)\sum_{i,j=1}^N z_iz_j\xi_i\xi_j \biggr) dt\\
&\ge \int_0^1 b(|y+t(x-y)|^2) |x-y|^2\,dt
\end{aligned}\end{equation}
where we have used \eqref{b2} and
$$\sum_{i,j=1}^N z_iz_j\xi_i\xi_j=\biggl(\sum_{i=1}^N z_i\xi_i\biggr)^2\ge 0.$$
Finally from \eqref{dif}, by using \eqref{xyt}, the fact that $t\mapsto b(t^2)$ is increasing and \eqref{growthb1} we get the claim, namely
$$\langle b(|x|^2)x-b(|y|^2)y, x-y \rangle\ge b\Big(\frac{1}{16}|x-y|^2\Big)|x-y|^2\ge C[|x-y|^2+|x-y|^q],$$
where $C>0$.
\end{proof}


\section{Problem \eqref{Pb} in the $L^2$-subcritical case}\label{app}

In this section, assuming that $N\ge 3$, $\rho>0$, $q>1$ and \eqref{h1.1}, we will prove Theorem \ref{mainb}, namely that \eqref{Pb} has a normalized ground state solution for a general function $b$ satisfying \eqref{b0}--\eqref{b2}, obtaining a global minimizer for 
$$ m(\rho)\!:=\!\inf _{u \in S_\rho} I(u)$$
where 
\begin{equation} \label{energyFb}
I(u)= \frac{1}{2}\int_{\mathbb R^N} B(|\nabla u|^2) \, dx-\frac{1}{p}{\|u \|}_p^p.
\end{equation}
Solutions of this type can be obtained by searching critical points of $I$
so that $\lambda$ in \eqref{Pb} appears as a Lagrange multiplier. This implies that $\lambda$ is part of the unknown.
Using \eqref{growthb2}, we can estimate the functional $I$ as follows
\begin{equation} \label{energyFinb}
\frac{c_1}{2} \bigl(\|\nabla u\|_2^2+\|\nabla u\|_q^q\bigr)-\frac{1}{p}{\|u \|}_p^p \le I(u)\le \frac{c_2}{2} \bigl(\|\nabla u\|_2^2+\|\nabla u\|_q^q\bigr)-\frac{1}{p}{\|u \|}_p^p
\end{equation}

Of course, the functional $I$ in \eqref{energyFb} is  well defined in the entire $\X$, since $2<p<\min\{2^*,q^*\}$.
The proof of the $\cC^1$ regularity of $I$ in $\X$ is almost standard.
In turn, $I': \X\to \X'$ is given by
\begin{equation}\label{I'b}
I'(u)[\phi]=\int_{\mathbb{R}^N} b(|\nabla u|^2)\nabla u \nabla \phi \, dx -\int_{\mathbb{R}^N}|u|^{p-2}u\phi \, dx
\end{equation}
for all $u, \phi\in \X$.
Note that \eqref{Pb} has always the trivial solution $0$, but the constraint $S_\rho$ prevents the case to occur.


 We start with the following lemma which justifies the global minimization problem.

\begin{lemma} \label{2.1.} For any $\rho>0$, then
$$-\infty< m(\rho)=\inf _{u \in S_\rho} I(u)<0.$$
\end{lemma}

\begin{proof}
Recalling that,  $\delta_r:=\frac{N(r-2)}{2r}$, for $r>1$,  observe that
\begin{equation*}
\left\{\begin{array}{ll}
p\delta_{p}<q(1+\delta_{q})<2&\mbox{if}~~~~\!q\!<\!2\!<\!p\!<\!q(1+\frac{2}{N}), 
\\
p\delta_{p}<2<q(1+\delta_{q})&\mbox{if}~~~~q\!>\!2~~~~\mbox{and}
~~~~2\!<\!p\!<\!2\!+\!\frac{4}{N},
\end{array}\right.
\end{equation*}
so we have
$$
   p\delta_{p}<\min\{ 2, q(1+\delta_{q}) \}.
$$
For any fixed $u \in S_{\rho}$, using \eqref{energyFinb} and the Gagliardo-Nirenberg inequality \eqref{equ2.2} we get
\begin{align} \label{lower-bounded1.1}
I(u)\ge\frac{c_1}{2} (\|\nabla u\|_{2}^2+\|\nabla u\|_{q}^q)-\frac{1}{p} {\|u\|}_p^p
\geq \frac{c_1}{2} \|\nabla u\|_{2}^2-\frac{\mathcal{C}^p_{N,p}}{p} \rho^{p(1-\delta_p)}\left\|\nabla u\right\|_{2}^{p\delta_p}.
\end{align} 
From \eqref{lower-bounded1.1}, we obtain $m(\rho)>-\infty$. Next, we show that $m(\rho)\!<\!0$. Let $u \in S_{\rho}$ be fixed, then we have $u_{t}(x)\!=\!t^{\frac{N}{2}}u(tx)\in S_{\rho}$ and
$$I(u_t)\le \frac{c_2 t^{2}}{2} \|\nabla u\|_{2}^2+\frac{c_2 t^{q(1+\delta_{q})}}{2} \|\nabla u\|_{q}^q-\frac{t^{ p\delta_{p}}}{p} {\|u\|}_p^p<0$$
for $t\!>\!0$ sufficiently small. This implies that $m(\rho)\!\leq\!I(u_t)\!<\!0$.
\end{proof}

In what follows we prove the continuity of $m(\rho)$ and the subadditivity property which will be useful in proving compactness.

\begin{lemma}\label{LemA3.5}
The following properties hold
\begin{enumerate}[label=(\roman{*}),ref=\roman{*}]
\item\label{i} the map $\rho\mapsto m(\rho)$ is continuous;
\item\label{ii} if $\hat \rho_1\!\in\!(0,\rho)$ and $\hat \rho_2\!=\!\sqrt{\rho^2-\hat \rho^2_1}$, we have $ m(\rho)< m(\hat \rho_1)+m(\hat \rho_2)$.
\end{enumerate}
\end{lemma}

\begin{proof}
\eqref{i}~~Let $\rho>0$ and $\{\rho_n\}\subset(0,+\infty)$ such that $\rho_n \to \rho$, it is sufficient to prove that $m(\rho_n) \to m(\rho)$. For every $n\in \mathbb{N}^+$, there exists $u_n \in S_{\rho_n}$ such that $m(\rho_n) \leq I(u_n) <m(\rho_n)+1/n$, where $m(\rho_n)\leq0$ from Lemma \ref{2.1.}. We first show that $\{u_{n}\}$ is bounded in $\X$. Similarly to \eqref{lower-bounded1.1}, we get
\begin{equation} \label{3.2boud2}
0\geq I(u_n)\ge\frac{c_1}{2} (\|\nabla u_n\|_{2}^2+ \|\nabla u_n\|_{q}^q)-\frac{1}{p} {\|u_n\|}_p^p \geq \frac{c_1}{2} \|\nabla u_n\|_{2}^2-\frac{\mathcal{C}^p_{N,p}}{p} \rho_n^{p(1-\delta_p)}\left\|\nabla u_n\right\|_{2}^{p\delta_p}.
\end{equation}
which gives
\begin{equation}\label{3.2bbb}
\|\nabla u_n\|_{2}\leq \biggl[\frac{2\mathcal{C}^p_{N,p}}{c_1 p} \biggr]^{\frac{1}{2-p\delta_p}}\rho_n^{\frac{p(1-\delta_p)}{2-p\delta_p}}\leq \biggl[\frac{2\mathcal{C}^p_{N,p}}{c_1 p} \biggr]^{\frac{1}{2-p\delta_p}}\rho^{\frac{p(1-\delta_p)}{2-p\delta_p}}+ o(1).
\end{equation}
With this upper bound, we derive from \eqref{3.2boud2} and the Gagliardo-Nirenberg inequality \eqref{equ2.2} that
\begin{equation}\label{3.2bi}\begin{aligned} \frac{c_1}{2} \|\nabla u_n\|_{q}^q\leq\frac{1}{p} {\|u_n\|}_p^p
&\leq \frac{\mathcal{C}^p_{N,p}}{p} \rho_n^{p(1-\delta_p)}\left\|\nabla u_n\right\|_{2}^{p\delta_p}\\&\leq \biggl(\frac{2}{c_1}\biggr)^{\frac{p\delta_p}{2-p\delta_p}}\biggl[\frac{\mathcal{C}^p_{N,p}}{p} \biggr]^{\frac{2}{2-p\delta_p}}\rho^{\frac{2p(1-\delta_p)}{2-p\delta_p}}+ o(1).
\end{aligned}\end{equation}
Therefore, $\{u_{n}\}$ is bounded in $\X$. Now considering $v_n := \frac{\rho}{\rho_n}u_n \in S_\rho$, we have
$$
m(\rho) \leq I(v_n)=\frac{\rho^p}{\rho_n^p}I(u_n)+\frac{1}{2}\int_{\mathbb R^N} \biggl[B\biggl(\frac{\rho^2}{\rho_n^2}|\nabla u_n|^2\biggr)-\frac{\rho^p}{\rho_n^p}B(|\nabla u_n|^2)\biggr]dx=I(u_n)+ o(1),
$$
where we used the boundedness of $\{u_n\}$, the continuity of $B$ and the fact that $\rho_n\to \rho$. Passing to the limit as $n \to +\infty$,
we deduce that
 $$  m(\rho) \leq \liminf_{n \to +\infty}  m(\rho_n).  $$
In a similar way, let $\{w_n\}$ be a minimizing sequence for $m(\rho)$, which is also bounded, and
let $z_n := \frac{\rho_n}{\rho}w_n \in S_{\rho_n}$. Then we have
$$
m(\rho_n) \leq I(z_n) = I(w_n) +  o(1) \Longrightarrow \limsup_{n \to +\infty} m(\rho_n) \leq m(\rho).$$
\eqref{ii}~~For any fixed $\hat \rho_1\in(0,\rho)$, we first claim that
\begin{equation} \label{CoTq4}
 m(\alpha \hat \rho_1)< \alpha^2 m(\hat \rho_1),\qquad\forall \alpha>1.
\end{equation}
Let $\{u_n\}\subset S_{\hat \rho_1}$ be a minimizing sequence for $m(\hat \rho_1)$, then $ u_n(\alpha^{-\frac{2}{N}}x) \in S_{\alpha \hat \rho_1}$. 
Take
\begin{align*}
m(\alpha \hat \rho_1)-\alpha^2I(u_n)&\leq I(u_n(\alpha^{-\frac{2}{N}}x))-\alpha^2I(u_n)\\
&=\frac{\alpha^2}{2}\int_{\mathbb R^N} \big(B(\alpha^{-4/N}|\nabla u_n(x)|^2)-B(|\nabla u_n(x)|^2)\big)dx.
\end{align*}
For all $\alpha>1$, then $\alpha^{-4/N}<1$ so that, by the monotonicity of $B$, we get
$$B(\alpha^{-4/N}|\nabla u_n(x)|^2)\le B(|\nabla u_n(x)|^2).$$
\\
As a consequence $m(\alpha \hat \rho_1)\leq\alpha^2 m(\hat \rho_1)$, with equality if and only if $\|\nabla u_n\|_{2}^2\to 0$ and $\|\nabla u_n\|_{q}^q\to0$ as $n\to+\infty$. In view of these facts, inequality \eqref{equ2.2} indicates that $\|u_n\|_{p}^p\to0$. It must be that $m(\alpha \hat \rho_1)<\alpha^2 m(\hat \rho_1)$, otherwise, we obtain the following contradiction
$$0 > m(\hat \rho_1) = \lim_{n\to+\infty} I(u_n)\ge  \frac{c_1}{2} \lim_{n\to+\infty}\|\nabla u_n\|_{2}^2+\frac{c_1}{2} \lim_{n\to+\infty}\|\nabla u_n\|_{q}^q-\frac{1}{p} \lim_{n\to+\infty}{\|u_n\|}_p^p= 0.$$
In the same manner, we can get
\begin{equation} \label{CoTq5}
 m(\alpha \hat \rho_2)< \alpha^2 m(\hat \rho_2),\qquad\forall \alpha>1.
\end{equation}
Finally, apply \eqref{CoTq4} with $\alpha=\dfrac{\rho}{\hat \rho_1}>1$ and \eqref{CoTq5} with $\alpha=\dfrac{\rho}{\hat \rho_2}>1$ respectively, 
we get
\begin{align*}
 m(\rho)\!=\!\frac{\hat \rho^2_1}{\rho^2}m\biggl(\frac{\rho}{\hat \rho_1}\hat c_1\biggr) \!+\!\frac{\hat \rho_2^2}{\rho^2}  m\biggl(\frac{\rho}{\hat \rho_2}\hat \rho_2\biggr)
\!< \!m(\hat \rho_1)\!+\!m(\hat \rho_2).
\end{align*}
\end{proof}

Applying Lemma \ref{LemA3.5}, we can prove the compactness of the minimizing sequences for $m(\rho)$.

\begin{lemma} \label{LeMa3.4.1}
Let $\{w_{n}\} \subset S_\rho$ be a minimizing sequence for $m(\rho)$,
then $m(\rho)$ possesses another minimizing sequence $\{u_{n}\} \subset S_\rho$ such that
$$\left\|u_{n}-w_{n}\right\|_\X \to 0,~~~~~~~~\big(I|_{S_\rho}\big)'\left(u_{n}\right) \to 0\qquad \text{as}~~~~n\to+\infty.$$
Moreover, $\{u_{n}\}$ is relatively compact in $\X$ up to translations and hence $m(\rho)$ is attained.
\end{lemma}


\begin{proof}
Since $\{w_{n}\} \subset S_\rho$ is a minimizing sequence of $m(\rho)$, from the Ekeland's variational principle (cf. \cite[Theorem2.4]{Mwlm}), we get a new minimizing sequence $\{u_{n}\} \subset S_\rho$ for $m(\rho)$ such that $\left\|u_{n}-w_{n}\right\|_{\X} \to 0$, which is also a Palais-Smale
sequence for $I|_{S_\rho}$. Hence, we have $\big(I|_{S_\rho}\big)'\left(u_{n}\right) \to 0$.
\\
In the same way as the proof of Lemma \ref{LemA3.5}-\eqref{i}, we obtain that $\{u_{n}\}$ is bounded in $\X$. If $\mathop {\lim }\limits_{n  \to +\infty} \sup _{y \in \mathbb{R}^{N}} \int_{B_{R}(y)}\left|u_{n}(x)\right|^{2}dx\!=\!0$ for any $R\!>\!0$, we can prove that $\left\|u_{n}\right\|_{r} \to 0$ for $2<r<2^{*}$ (see \cite[Lemma I.1]{LilP}). In particular, we have ${\|u_n\|}_p^p\to0$, this together with Lemma \ref{2.1.} leads to the following contradiction
\begin{align*}
0&>m(\rho)=\mathop {\lim }\limits_{n  \to +\infty}I(u_n)\ge\frac{c_1}{2}\mathop {\lim }\limits_{n  \to +\infty}{\|\nabla u_n\|}_2^2+\frac{c_1}{2} \mathop {\lim }\limits_{n  \to +\infty}\|\nabla u_n\|_{q}^q-\frac{1}{p}\mathop {\lim }\limits_{n  \to +\infty}{\|u_n\|}_p^p \geq0.
\end{align*}
Then, there exist an $\varepsilon_0>0$ and a sequence $\{y_{n}\} \subset \mathbb{R}^{N}$ such that
$$
\int_{B_{R}(y_n)}\left|u_{n}(x)\right|^{2} d x \geq \varepsilon_0>0
$$
for some $R>0$. Hence we have $u_{n}(x+y_n)\rightharpoonup u_\rho\not\equiv0$ in $\X$ for some $u_\rho\in \X$. Since the problem is invariant by translation, without loss of generality we can assume that $y_n=0$. Let $v_n:=u_n -u_\rho$, then we see that $v_{n}\rightharpoonup 0$ in $\X$. Therefore, we get
\begin{equation}
\begin{aligned} \label{NlaG1} &\left\|u_{n}\right\|_{2}^2
=\left\|v_{n}\right\|_{2}^2
+\left\|u_\rho\right\|_{2}^2+ o(1), \\
&\left\|\nabla u_{n}\right\|_{2}^2
=\left\|\nabla v_{n}\right\|_{2}^2
+\left\|\nabla u_\rho\right\|_{2}^2+ o(1).
\end{aligned}
\end{equation}
Moreover, by the compactness of the embedding of $D^{1,q}(\mathbb R^N)$ in $L^s_{loc}(\mathbb R^N)$ for any $1\le s<q^*$ if $1<q<N$ or $1\le s<+\infty$ if $q\ge N$, we have
$v_n\to 0$ in $L^{s}(\omega)$ with $\omega\Subset\mathbb{R}^N$. So, by using an increasing sequence of compact sets whose union is $\mathbb R^N$ and a diagonal argument, we get $v_n(x)\to 0$ for a.e. $x\in\mathbb R^N$.
From the Br\'{e}zis-Lieb Lemma \cite{a7}, we have
\begin{equation}\label{plieb}
\left\|u_{n}\right\|_{p}^p
=\left\|v_{n}\right\|_{p}^p
+\left\|u_\rho\right\|_{p}^p+ o(1).
\end{equation}
By a standard truncation argument (cf. \cite{aLyM,ByF2,Wlzk}) in what follows we prove that,  up to subsequences, 
\begin{equation}\label{aeconv}
\nabla  u_{n} \rightarrow \nabla u_\rho \quad \text{a.e. on} \,\,\R^N.
\end{equation}
Choose $\psi\in C_0^\infty(\mathbb R^N)$ such that
$0\le\psi\le 1$ in $\mathbb R^N$, $\psi(x)=1$ for every $x\in B_1(0)$ and $\psi(x)=0$ for every $x\in \mathbb R^N\setminus B_2(0)$. Now, take $R>0$ and define $\psi_R(x)=\psi(x/R)$ for $x\in\mathbb R^N$.
Since $\{u_n\}\subset S_\rho$ and $\big(I|_{S_\rho}\big)'\left(u_{n}\right) \to 0$, by \cite[Lemma 3]{BerLionsII}, we have that
\begin{equation*}
I'(u_n)-\frac{1}{\rho^2}I'(u_n)[u_n]u_n \to 0, \qquad\text{in }\X',
\end{equation*}
and so
\begin{equation*}
I'(u_n)[(u_n-u_\rho)\psi_R]-\frac{1}{\rho^2}I'(u_n)[u_n]\langle u_n, (u_n-u_\rho)\psi_R\rangle_{L^2} \to 0.
\end{equation*}
Since $\{u_n\}$ is bounded and  $\psi_R$ has compact support,
we have $I'(u_n)[(u_n-u_\rho)\psi_R]\to0$ as $n\to+\infty$.
Recalling the definition of $I'(u)[\phi]$ in \eqref{I'b} with $u=u_n$ and $\phi=(u_n-u_\rho)\psi_R$, we get
\begin{equation}\label{qog}\begin{aligned}
\int_{\mathbb R^N} &\biggl[b(|\nabla u_n|^2)\nabla u_n-b(|\nabla u_\rho|^2)\nabla u_\rho\biggr](\nabla u_n-\nabla u_\rho)\psi_R \, dx=I'(u_n)[(u_n-u_\rho)\psi_R]\\
&-\int_{\mathbb R^N}b(|\nabla u_n|^2)\nabla u_n u_n  \nabla \psi_R \, dx+\int_{\mathbb R^N}b(|\nabla u_n|^2)\nabla u_n u_\rho  \nabla \psi_R \, dx\\
&+ \int_{\mathbb R^N} |u_n|^p \psi_R\,  dx-\int_{\mathbb R^N} |u_n|^{p-2}u_n u_\rho \psi_R \, dx\\
&-\int_{\mathbb R^N}b(|\nabla u_\rho|^2)\nabla u_\rho(\nabla u_n-\nabla u_\rho)\psi_R\,  dx.
\end{aligned}\end{equation}
Moreover, by the definition of $\psi_R$ and the boundedness of $u_n, u_\rho\in \X$, then
$$\int_{\mathbb R^N}\nabla u_n 
u_n\nabla \psi_R  \, dx\to \int_{\mathbb R^N} \nabla u_n u_\rho \nabla \psi_R\, dx,$$
$$\int_{\mathbb R^N} |\nabla u_n|^{q-2}\nabla u_n u_n \nabla \psi_R\, dx\to \int_{\mathbb R^N} |\nabla u_n|^{q-2}\nabla u_n u_\rho \nabla \psi_R\,  dx,$$
which imply
$$\int_{\mathbb R^N}b(|\nabla u_n|^2)\nabla u_n u_n  \nabla \psi_R \, dx \to\int_{\mathbb R^N}b(|\nabla u_n|^2)\nabla u_n u_\rho  \nabla \psi_R \, dx, \qquad n\to+\infty.$$
In view of the weak convergence of $u_n \weakto u_\rho$ in $\X$,
it holds
$$\int_{\mathbb R^N} |\nabla u_\rho|^{q-2}\nabla u_\rho \nabla u_n \psi_R \, dx\to \int_{\mathbb R^N} |\nabla u_\rho|^{q} \psi_R \, dx$$
$$\int_{\mathbb R^N} \nabla u_n \nabla u_\rho \psi_R \, dx\to \int_{\mathbb R^N} |\nabla u_\rho|^2 \psi_R \, dx,$$
which implies
$$\int_{\mathbb R^N}b(|\nabla u_\rho|^2)\nabla u_\rho(\nabla u_n-\nabla u_\rho)\psi_R \, dx\to 0, \qquad n\to+\infty.$$
Moreover
$$\int_{\mathbb R^N} |u_n|^p \psi_R \, dx, \int_{\mathbb R^N} |u_n|^{p-2}u_n u_\rho \psi_R\,  dx\to \int_{\mathbb R^N} |u_\rho|^{p}\psi_R \, dx, \qquad n\to+\infty.$$
So that, from \eqref{qog} by letting $n\to+\infty$, we proved
$$\lim_{n\to+\infty}\int_{\mathbb R^N} \biggl[b(|\nabla u_n|^2)\nabla u_n-b(|\nabla u_\rho|^2)\nabla u_\rho\biggr](\nabla u_n-\nabla u_\rho)\psi_R \, dx=0.$$
By virtue of Lemma \ref{lemdis}, we get
$$\nabla u_n\to\nabla u_\rho \quad \text{in} \,\,\, L^2(B_{2R})\cap L^q(B_{2R}).$$
Finally, by the arbitrariness of $R$, then \eqref{aeconv} is reached. 
\\
We next claim that $u_{n}\rightarrow u_\rho\not\equiv0~~~~\mbox{in}~~~~L^2(\R^N)$, or equivalently $v_{n}\!\rightarrow \!0$ in $L^2(\R^N)$. Denote $\left\|u_\rho\right\|_{2}\!=\!\hat \rho_1$. If $\hat \rho_1\!=\!\rho$, the proof is completed by \eqref{NlaG1}, now assume $\hat \rho_1\!<\!\rho$ and so, by \eqref{NlaG1}, $\|v_n\|_2 \to \hat \rho_2$, where  $\hat \rho_2\!=\!\sqrt{\rho^2-\hat \rho^2_1}>0$.
Now, applying Brezis-Lieb Lemma, namely Theorem 2 in \cite{a7} with $j(s)=B(s^2)$ 
which is convex by \eqref{b2}
(see Example (b) in \cite{a7}) by \eqref{aeconv}, and using \eqref{plieb}, we have
$$m(\rho)=I\left(u_{n}\right)+ o(1)
=I\left(v_{n} \right)+I\left(u_\rho \right)+ o(1)\geq m(\left\| v_{n}\right\|_{2})+m(\hat \rho_1).$$
By the continuity of $\rho\mapsto m(\rho)$ (see Lemma \ref{LemA3.5} (i)), we have
\begin{align} \label{NlaG4}
m(\rho)\geq m(\hat \rho_2)+m(\hat \rho_1).
\end{align}
However, taking 
$\hat \rho_1<\rho$, \eqref{NlaG4} contradicts to  Lemma \ref{LemA3.5} (ii). 
Therefore, we have $\left\|u_\rho\right\|_{2}\!=\!\rho$ and hence $v_{n}\rightarrow 0$ in $L^2(\R^N)$. It follows immediately, by \eqref{equ2.2}, that
$$\|v_n\|_{p} \leq \mathcal{C}_{N,p} \left\|\nabla v_n\right\|_{2}^{\delta_p} \left\|v_n\right\|_{2}^{(1-\delta_p)} \to0.$$
Finally, we get 
\begin{align*}
m(\rho)=I\left(v_{n} \right)+I\left(u_\rho \right)+ o(1)\geq\frac{c_1}{2}{\|\nabla v_{n}\|}_2^2+\frac{c_1}{2} \|\nabla v_n\|_{q}^q+m(\rho)+ o(1),
\end{align*}
which indicates ${\|\nabla v_{n}\|}_2=  o(1)$ and ${\|\nabla v_{n}\|}_q=  o(1)$. So we have $v_{n}\!\rightarrow \!0$ in $\X$ and $u_{n}\rightarrow u_\rho\not\equiv0$ in $\X$.
\end{proof}

Now we are ready to prove Theorem \ref{mainb}, whose statement is given in the Introduction.

\begin{proof}[Proof of Theorem  \ref{mainb}.]
By using Lemma \ref{LeMa3.4.1}, we see that
$m(\rho)$ is attained by some $u_\rho \!\in\! S_\rho $. Next, the Lagrange multipliers rule implies the existence of some Lagrange multiplier $\lambda_\rho\in \mathbb{R}$ such that  $(\lambda_\rho,u_\rho)$ satisfies \eqref{Pb}. From \cite [Theorem 2.7]{AL}, we have $I\left({|u_\rho|}^*\right)\!\leq\!I\left(|u_\rho|\right)\!\leq\!I\left(u_\rho\right)$, where ${|u_\rho|}^*$ is the symmetric decreasing rearrangement of $|u_\rho|$. So we can assume that $u_\rho\!\in\! S_\rho$ is nonnegative and radially decreasing. Thus, $u_\rho$ is a ground state solution also among radial solutions and Theorem \ref{mainb} is proved. 
\end{proof}



\section{Born-Infeld problem \eqref{P} in the  $L^2$-subcritical case}\label{res}

This section is devoted to get a solution to the Born-Infeld problem \eqref{P} where $a$ satisfies \eqref{a0}--\eqref{a2}. First, we have to truncate the map $a$, as in a similar way of \cite{BDD,MP}, to obtain an operator which behaves like $b$ in \eqref{Pb}. Then by applying Theorem \ref{mainb} to the truncated problem, we can find a normalized solution for such a problem. Finally, in order to get rid of the truncation we have to restrict ourself to the radial setting and, moreover, uniform bounds of the solution to the truncated problem are required.

Let us fix $\t_1\in(0, 1)$ and, for any $\t\in(0,\t_1]$, we take $q=q(\t)>N$. 
Define $a_\t:[0,+\infty)\to \R^+$ a $\cC^1$ function defined as follows:
$$a_\t(s)\!:=\!
\begin{cases}
a(s),& s\le 1-\t,
\\
(1-\theta)^{-\frac{q-4}{2}}\dfrac{2}{q-2} a'(1-\theta)\big(s^\frac{q-2}{2}-(1-\theta)^\frac{q-2}{2}\big)+a(1-\theta), & s> 1-\t,
\end{cases}$$
and set $A_\t(s):=\int_0^sa_\t(t)\, dt$.
Observe that  $a_\t, A_\t$ satisfy the following growth conditions
$$c_1(s^2+|s|^q)\le a_\t(s^2)s^2\le c_2(s^2+|s|^q)$$
\begin{equation}\label{growth2}
c_1(s^2+|s|^q)\le A_\t(s^2)\le c_2(s^2+|s|^q)
\end{equation}
where $c_1>0$ does not depend on $\t$ and $c_2=c_2(\t)>0$. 
Then, we consider the truncated problem 
\begin{equation}\label{Ptr}\tag{${\mathcal P}_\t$}
\begin{cases}
-\mbox{div}(a_\t(|\nabla u|^2)\nabla u)=|u|^{p-2}u-\lambda u, \qquad & \text{in}\,\,\mathbb R^N,\\
u\in \mathcal X_r,
\\
\int_{\R^N}|u|^2\,dx=\rho^2.
\end{cases}\end{equation}
Clearly, if $u_\t$ is a solution of \eqref{Ptr} such that $|\nabla u_\t|^2 \le 1 - \t$, then $u_\t$ is a solution also of \eqref{P}.
\\

We are looking for normalized solutions to \eqref{Ptr} under the constraint $S_\rho$ for fixed $\rho>0$, defined in \eqref{Sc}.
Solutions of this type can be obtained by searching critical points of the following functional
\begin{equation} \label{energyF}
I_\t(u)= \frac{1}{2}\int_{\mathbb R^N} A_\t(|\nabla u|^2) dx-\frac{1}{p}{\|u \|}_p^p
\end{equation}
so that $\lambda$ in \eqref{Ptr} appears as a Lagrange multiplier. Moreover, we define
$$ m_\t(\rho):=\inf _{u \in S_\rho} I_\t(u).$$
This implies that $\lambda$ is part of the unknown. Similarly for $I$ in Section \ref{app}, using \eqref{growth2} we can estimate the functional $I_\t$ as follows
\begin{equation} \label{energyFin}
\frac{c_1}{2} \bigl(\|\nabla u\|_2^2+\|\nabla u\|_q^q\bigr)-\frac{1}{p}{\|u \|}_p^p \le I_\t(u)\le \frac{c_2}{2} \bigl(\|\nabla u\|_2^2+\|\nabla u\|_q^q\bigr)-\frac{1}{p}{\|u \|}_p^p
\end{equation}

Of course, the functional $I_\t$ in \eqref{energyF} is  well defined in the entire $\X_r$, since  $2<p<\min\{2^*,q^*\}=2^*$.
The proof of the $\C^1$ regularity of $I_\t$ in $\X_r$ is almost standard.
In turn, $I_\t': \X_r\to (\X_r)'$ is given by
$$I_\t'(u)[\phi]=\int_{\mathbb{R}^N} a_\t(|\nabla u|^2)\nabla u \nabla \phi dx -\int_{\mathbb{R}^N}|u|^{p-2}u\phi dx$$
for all $u, \phi\in \X_r$.
Note that \eqref{Ptr} has always the trivial solution $0$, but the constraint $S_\rho$ prevents the case from occurring.

The following is a crucial step in our arguments.

\begin{proposition}\label{pr:solz}
	For any $\t\in (0,\t_1]$, there exists $u_\t\in \X_r$ a non-trivial normalized  solution of \eqref{Ptr}, for a certain $\lambda_\t\in \R$, such $I_\t(u_\t)=m_\t(\rho)$. Moreover
	there exists $C>0$ such that 
	\begin{equation}\label{unifbddz}
	\|u_\t\|_{\X}\le C\text{ and }|\lambda_\t|\le C, \quad\hbox{ for any }\t\in (0,\t_1].
	\end{equation}
	Finally $u_\t$ is a weak solution of 
	\begin{equation}\label{eqradz}
	-\big(r^{N-1}a_\t(|u'_\t(r)|^2)u'_\t(r)\big)'=r^{N-1}[|u_\t(r)|^{p-2}u_\t(r)-\lambda_\t u_\t(r)],
	\end{equation}
	namely 
	\[
	\int_0^{+\infty}r^{N-1}a_\t(|u'_\t(r)|^2)u'_\t(r)v'(r)\, dr
	=\int_0^{+\infty}r^{N-1}[|u_\t(r)|^{p-2}u_\t(r)-\lambda_\t u_\t(r)]v(r)\,dr,
	\]
	for all $v\in \X_r$.
\end{proposition}

\begin{proof}
Observe that $a_\t$ satisfies \eqref{b0}--\eqref{b2} and so the results of
Theorem \ref{mainb} hold also for \eqref{Ptr}. So there exists $u_\t\in \X_r$ nonnegative and radially decreasing such that
\[
I_\t(u_\t)=m_\t(\rho)=\inf _{u \in S_\rho} I_\t(u)<0.
\] 
 It remains to prove \eqref{unifbddz}.
Using \eqref{energyFin}, the Gagliardo-Nirenberg inequality \eqref{equ2.2}  and Lemma \ref{2.1.} we get
$$0>m_\t(\rho)=I_\t(u_\t)\ge\frac{c_1}{2} (\|\nabla u_\t\|_{2}^2+\|\nabla u_\t\|_{q}^q)-\frac{1}{p} {\|u_\t\|}_p^p.$$
Following the same ideas to get \eqref{3.2bbb} and \eqref{3.2bi} we arrive to
$$\|\nabla u_\t\|_{2}\leq \biggl[\frac{2\mathcal{C}^p_{N,p}}{c_1 p} \biggr]^{\frac{1}{2-p\delta_p}}\rho^{\frac{p(1-\delta_p)}{2-p\delta_p}},\qquad
\|\nabla u_\t\|_{q}^q\leq \biggl(\frac{2}{c_1}\biggr)^{\frac{2}{2-p\delta_p}}\biggl[\frac{\mathcal{C}^p_{N,p}}{p} \biggr]^{\frac{2}{2-p\delta_p}}\rho^{\frac{2p(1-\delta_p)}{2-p\delta_p}},$$
which gives the first part of \eqref{unifbddz}, since $u_\t\in S_\rho$. Moreover, multiplying the equation of \eqref{Ptr} by $u_\t$ and integrating, since we have already proved that $\{u_\t\}$ is uniformly bounded, with respect to $\t$, we obtain also the uniform bound on $\{\lambda_\t\}$.
\end{proof}

We are now able to conclude the proof of Theorem \ref{main}.

\begin{proof}[Proof of Theorem \ref{main}]
For any $\t\in (0,\t_1]$, there exists $u_\t\in \X_r$ a non-trivial nonnegative and radially decreasing normalized  solution of \eqref{Ptr}, for a certain $\lambda_\t\in \R$, such $I_\t(u_\t)=m_\t(\rho)$. Since $q>N$, $u_\t\in L^\infty(\R^N)$ and $u_\t$ is a solution of \eqref{eqradz} in $(0,+\infty)$, it is easy to check that $u_\t$ is regular for $r>0$. 
Moreover,  $r^{N-1}a_\t(|u'_\t(r)|^2)u'_\t(r)$ satisfies the Cauchy condition at the origin so that it has a finite limit as $r\to0$. 
Let us prove the following.
\\
{\sc Claim $1$:} we have that
\begin{equation}\label{eq:lim0}
\lim_{r\to 0}r^{N-1}a_\t(|u'_\t(r)|^2)u'_\t(r)=0.
\end{equation}
Suppose, by contradiction, that it is different from zero, then $\lim_{r\to0}|u'_\t(r)|=+\infty$. By using \eqref{growthb1} and since $u_\t$ is constant on a sphere centered at $0$, in view of Lieberman's result \cite{Lieb}, then $u_\t\in \cC^{1,\alpha}$ in a neighborhood of $0$ for some $\alpha\in(0,1)$, cf. \cite{MP}. This is in contradiction with $\lim_{r\to0}|u'_\t(r)|=+\infty$.
\\
{\sc Claim $2$:} there exists $C>0$ such that 
\begin{equation}\label{claim1}
|a_\t(|u'_\t(r)|^2)u'_\t(r)|\le C, \qquad \hbox{for any $r\ge 0$ and $\t\in (0,\t_1]$}.
\end{equation}
Integrating the equation \eqref{eqradz} for any $r_2>r_1>0$, we have
\begin{equation*}
-a_\t(|u'_\t(r_2)|^2)u'_\t(r_2)+\frac{r_1^{N-1}}{r_2^{N-1}}a_\t(|u'_\t(r_1)|^2)u'_\t(r_1)
=\frac{1}{r_2^{N-1}}\int_{r_1}^{r_2} s^{N-1}\left(|u_\t(s)|^{p-2}-\lambda_\t\right)u_\t(s)\, ds,
\end{equation*}
and, letting $r_1\to 0$, by \eqref{eq:lim0}, we deduce that
\begin{equation*}
\left|a_\t(|u'_\t(r_2)|^2)u'_\t(r_2)\right|\le
\frac{1}{r_2^{N-1}}\int_{0}^{r_2} s^{N-1}|u_\t(s)|^{p-1}\, ds
+\frac{|\lambda_\t|}{r_2^{N-1}}\int_{0}^{r_2} s^{N-1}|u_\t(s)|\, ds.
\end{equation*}
Observe that, by \eqref{unifbddz} and $\{u_\t\}$ is uniformly bounded in $L^\infty$, with respect to $\t\in (0,\t_1]$
$$\frac{1}{r_2^{N-1}}\int_{0}^{r_2} s^{N-1}|u_\t(s)|^{p-1}\, ds\le \frac{C}{N}r_2,\qquad \frac{1}{r_2^{N-1}}\int_{0}^{r_2} s^{N-1}|u_\t(s)|\, ds\le \frac{C}{N}r_2.$$
Therefore
\[
\lim_{r\to 0}a_\t(|u'_\t(r)|^2)u'_\t(r)=0,
\]
hence
\[
\lim_{r\to 0}u'_\t(r)=0.
\]
So, integrating the equation \eqref{eqradz}, for any $r>0$, we have
\begin{equation*}
-a_\t(|u'_\t(r)|^2)u'_\t(r)
=\frac1{r^{N-1}}\int_0^r s^{N-1}[|u_\t(s)|^{p-2}u_\t(s)-\lambda_\t u_\t(s)]\, ds.
\end{equation*}
By the continuous embedding of $\X$ in $L^s(\RN)$, for  $s\in [2,+\infty]$, and \eqref{unifbddz},  there exists $C>0$ such that $\|u_\t\|_s\le C\|u_\t\|_\X\le C$, for $s\in [2,+\infty]$ and any $\t\in (0,\t_1]$.  So, we have that, for any $0<r\le 1$ and $\t\in (0,\t_1]$,  
\begin{align*}
|a_\t(|u'_\t(r)|^2)u'_\t(r)|
&\le \frac{1}{r^{N-1}}\int_0^r s^{N-1}\big||u_\t(s)|^{p-2}u_\t(s)-\lambda_\t u_\t(s)\big|\, ds
\\&\le \frac{1}{r^{N-1}}\biggl[\int_0^r s^{N-1}|u_\t(s)|^{p-1}\, ds +|\lambda_\t|\int_0^r s^{N-1}| u_\t(s)|\, ds\biggr]
 \le Cr\le C.
	\end{align*}
By Proposition \ref{pr:solz}.
	If $r> 1$, then
	\begin{align*}
	|a&_\t(|u'_\t(r)|^2)u'_\t(r)|
	\le \frac1{r^{N-1}}\int_0^r s^{N-1}\big||u_\t(s)|^{p-2}u_\t(s)-\lambda_\t u_\t(s)\big|\, ds
	\\
	&\le \frac1{r^{N-1}}\biggl(\int_0^{1} s^{N-1}\big||u_\t(s)|^{p-2}u_\t(s)-\lambda_\t u_\t(s)\big|\, ds
	\\&\hspace{3cm}+\int_{1}^r s^{N-1}\big||u_\t(s)|^{p-2}u_\t(s)-\lambda_\t u_\t(s)\big|\, ds\biggr)\\
	&\le \frac C{r^{N-1}}
	+\underbrace{\frac{1}{r^{N-1}}\int_1^r s^{N-1}\big||u_\t(s)|^{p-2}u_\t(s)-\lambda_\t u_\t(s)\big|\, ds}_{(A)}.
	\end{align*}
	We have to estimate $(A)$. First of all, by Lemma \ref{strauss} and \eqref{unifbddz}, for $r>1$, we have that
	\[
	|u_\t(r)| \le C r^{-\frac{N-1}2}\|\n u_\t\|_2\le C r^{-\frac{N-1}2}.
	\]
	Hence, since $p>2$,
	\begin{align*}
	(A)&
	\le \frac{C}{r^{N-1}}\int_1^r s^{N-1}\big(|u_\t(s)|^{p-1}+|\lambda_\t||u_\t(s)|\big)\, ds
	\\
	&\le \frac{ C}{r^{N-1}}\int_1^r s^{N-1-\frac{N-1}2}\, ds
	\le C \left(r^{1-\frac{N-1}2}+1\right)\le C.
	\end{align*}
Since $N\ge 3$ the claim  is proved.
\\
In order to get a solution of the original problem \eqref{P}, it remains to prove the following.
\\
{\sc Claim $3$:} There exists $\bar\t\in (0, \t_1]$ such that
\begin{equation}\label{claim2}
|u'_\t(r)|\le 1-\bar\t, \quad \text{for any} \,\, r\ge0.
\end{equation}
Suppose by contradiction that \eqref{claim2} does not hold, then there exists a sequence $\{\t_n\}\subset(0, \t_1]$
which tends to zero and a sequence $\{r_n\}\subset \R^+$ such that
$$\lim_{n\to+\infty}|u'_{\t_n}(r_n)|=1,$$
which implies, by \eqref{a1}, that
$$\lim_{n\to+\infty}|a_{\t_n}(|u'_{\t_n}(r_n)|)u'_{\t_n}(r_n)|=+\infty.$$
Thus we obtain a contradiction with \eqref{claim1}.

\end{proof}

\begin{remark}






In general we cannot say that $u_\t$, the solution found in Theorem \ref{main}, is a ground state for \eqref{P} but,  whenever 
\begin{equation}\label{cla}
a_\t(s)\leq a(s), \qquad\text{ for }s\in (1-\theta,1),
\end{equation}
then we assert the following
$$I(u_\theta)=\inf_{S_\rho\cap \{\|\nabla u\|_\infty\leq 1\}}I.$$
Indeed, under the condition \eqref{cla}, we  have
$$\inf_{S_\rho\cap \{\|\nabla u\|_\infty\leq 1\}}I\leq I(u_\theta)=I_\theta(u_\theta)=\inf_{S_\rho}I_\theta\leq \inf_{S_\rho\cap \{\|\nabla u\|_\infty\leq 1\}}I_\theta\leq \inf_{S_\rho\cap \{\|\nabla u\|_\infty\leq 1\}}I.$$
Observe that \eqref{cla} holds if, for any $s\in (1-\t,1)$,  
$a'_\t(s)\le a'(s)$, which is equivalent to require that $ (1-\theta)^{-\frac{q-4}{2}} a'(1-\theta)\le s^{-\frac{q-4}{2}} a'(s)$, and so that
$ s^{-\frac{q-4}{2}} a'(s)$ is increasing. This is verified if \begin{equation*}
a''(s)s-\frac{q-4}{2}a'(s)>0,\qquad\text{ for }s\in (1-\theta,1).
\end{equation*}
If $q=4$, it holds if $a$ is assumed convex. Trivially, it is satisfied by the Born-Infeld operator \eqref{biop}.


\end{remark}

\section{Problem \eqref{Pb} in the  $L^2$-supercritical case}\label{super}

This section is devoted to the problem \eqref{Pb}, where $b$ satisfies \eqref{b0}--\eqref{b3}, in the  $L^2$-supercritical case. Indeed, in what follows we  assume \eqref{eq:L2super}.
Pursuing the approach in \cite{HIT}, we need to introduce an auxiliary functional $J$.

For any $u\in H^1(\mathbb R^N)$ and $\s\in\mathbb R$, we define 
$$(\s\ast u)(x):=e^{\sigma N/2}u(e^{\s} x)$$
and the  functional $J:\R\times \X_r\to \R$  as 
\begin{equation}\label{Jtheta}
J(\s, u)=
I\big(e^{\sigma N/2}u\big(e^{\s}\cdot )\big)
=\frac{e^{-N\s}}2 \irn B\big(e^{\s(N+2)}|\n u|^2\big)\, dx-\frac{e^{\s N(p/2-1)}}{p}\|u\|_p^p.
\end{equation}

Let us start with a characterization of the auxiliary functional $J$.
\begin{lemma}\label{lemA1}
If $u\neq 0$, then
$$J(\s, u)=I\big(e^{\sigma N/2}u(e^{\s}\cdot )\big)\to -\infty$$
as $\sigma\to +\infty$.
\end{lemma}

\begin{proof}
Note that, from \eqref{growthb2}, it holds
$$\begin{aligned}
c_1&\left(e^{\s(N+2)}|\n u|^2+e^{\s q(N+2)/2}|\n u|^q\right)\le B\big(e^{\s(N+2)}|\n u|^2\big)\\
&\hspace{2cm}\le c_2\left(e^{\s(N+2)}|\n u|^2+e^{\s q(N+2)/2}|\n u|^q\right),
\end{aligned}$$
so that
$$\begin{aligned}c_1&\left(e^{2\s}|\n u|^2+e^{\s q(N+2)/2-N\s}|\n u|^q\right)\le e^{-N\s}B\big(e^{\s(N+2)}|\n u|^2\big)\\
&\hspace{2cm}\le c_2\left(e^{2\s}|\n u|^2+e^{\s q(N+2)/2-N\s}|\n u|^q\right).\end{aligned}$$
Thus, taking into account \eqref{Jtheta} and \eqref{eq:L2super}, the claim is proved.
\end{proof}

%
%
%
%


Now we are ready to show that $I|_{S_{\rho,r}}$ satisfies the mountain pass geometry.

\begin{lemma}\label{IMP}
The functional $I|_{S_{\rho,r}}$ has a mountain pass geometry, namely
\begin{enumerate}[label=(\roman{*}),ref=\roman{*}]
\item\label{IMP1}  there exists $\eta>0$ such that
\begin{equation}\label{step1}
\inf_{V_\eta}I>0,
\end{equation}
where 
$$V_\eta:=\{u\in S_{\rho,r}: Q(u)=\eta\}\quad\text{and} \quad   Q(u)= \|\nabla u\|_2^2+\|\nabla u\|_q^q;$$
\item\label{IMP2} there exist $u_0,u_1\in S_{\rho,r}$ 
such that
$$Q(u_0)<\eta<Q(u_1), \quad I(u_0)<\inf_{V_\eta}I, \quad I(u_1)<0.$$
\end{enumerate}
\end{lemma}

\begin{proof}
\eqref{IMP1} Observe that, by using \eqref{growthb2} and \eqref{equ2.2},  for any $u\in S_{\rho,r}$ we have
$$\begin{aligned}
I(u)&= \frac{1}{2}\int_{\mathbb R^N} B(|\nabla u|^2) dx-\frac{1}{p}{\|u \|}_p^p\ge \frac{1}{2}c_1Q(u)
-\frac{1}{p}{\|u \|}_p^p\\
&\ge \frac{1}{2}c_1Q(u)-\frac{1}{p}\mathcal{C}_{N,p}^p \rho^{p(1-\delta_p)}\left\|\nabla u\right\|_{2}^{p\delta_p} \ge \frac{1}{2}c_1Q(u)-\frac{1}{p}\mathcal{C}_{N,p}^p \rho^{p(1-\delta_p)}Q(u)^{\frac{p\delta_p}2}. 
\end{aligned}$$
Since $p\delta_p>2$ by \eqref{h1.1}, we can take $\eta$ sufficiently small, namely $0<\eta<\bar\eta_\rho$ where
$$\bar\eta_\rho=\left[\frac{2}{pc_1}\mathcal{C}_{N,p}^p \rho^{p(1-\delta_p)}\right]^{\frac 2{2-p\delta_p}},$$
so that \eqref{step1} holds.
\\
\eqref{IMP2} Fix $u\in S_{\rho,r}$. Observe that
$$Q(\sigma\ast u)=e^{2\s}\|\nabla u\|_2^2+e^{\s\left(q\frac{N+2}{2}-N\right)}\|\nabla u\|_q^q,$$
and
$$ q\frac{N+2}{2}-N>0 \iff q>\frac{2N}{N+2},$$
which holds from \eqref{eq:L2super}. Now, by \eqref{step1} and $I(0)=0$ from the continuity of $I$, there is $\sigma_0<0$, with $|\sigma_0|$ sufficiently large, such that
$I(u_0)<\inf_{V_\eta}I$ and $Q(u_0)<\eta$, where $u_0:=\sigma_0\ast u\in S_{\rho,r}$.
\\
Moreover, there exists $\sigma_1>0$ sufficiently large such that, denoting $u_1:= \sigma_1\ast u\in S_{\rho,r}$ , then we have $I(u_1)<0$, by Lemma \ref{lemA1}, and 
$$Q(u_1)=e^{2\s}\|\nabla u\|_2^2+e^{\s\left(q\frac{N+2}{2}-N\right)}\|\nabla u\|_q^q>\eta.$$
\end{proof}

Since 
$I|_{S_{\rho,r}}$ has a mountain pass geometry, we can define its mountain pass level as
$$m:=\inf_{\g\in  \Gamma}\max_{t\in [0,1]}I\big(\g(t)\big),$$
where
\begin{equation}\label{gammamp}
\Gamma:=\{\g\in \cC([0,1],S_{\rho,r})\mid \g(0)=u_0\hbox{ and }\g(1)=u_1\}.
\end{equation}
Observe that 
\begin{equation}\label{muou1}
m\ge \inf_{V_\eta}I>\max\{I(u_0),I(u_1)\}.
\end{equation}
We define the following minimax level for $J|_{\R\times S_{\rho,r}}$:
\[
\tilde m:=\inf_{(\s,\g)\in \Sigma\times \Gamma}\max_{t\in [0,1]}J\big(\s(t),\g(t)\big), 
\]
where $\Gamma$ is defined in \eqref{gammamp} and 
$$\Sigma:=\{\s\in \cC([0,1],\R)\mid \s(0)=\s(1)=0\}.$$

Observe that it holds the following
\begin{lemma}\label{mm}
The minimax levels of $I$ and $J$ coincide, namely $m=\tilde m$.
\end{lemma}

\begin{proof}
For any $\gamma\in\Gamma$ we can see that $(0,\gamma)\in\Sigma\times\Gamma$ so, regarding $\Gamma\subset\Sigma\times\Gamma$, we have $m\ge\tilde m$.
Next, for any given $\big(\s(t),\gamma(t)\big)\in \Sigma\times\Gamma$, we set $\tilde\gamma(t)(x)=e^{\s(t)N/2}\gamma(t)(e^{\s(t)}x)$. We can verify that $\tilde\gamma\in\Gamma$ and, by \eqref{Jtheta}, that $I(\tilde\gamma(t))=J\big(\s(t),\gamma(t)\big)$. Thus we also have $m\le\tilde m$.
\end{proof}

Note that we are dealing with $J$ restricted to $\R\times S_{\rho,r}\subset \R\times \X_r$, so we recall that the norm of the derivative in the restriction is defined as
\begin{equation}\label{rema}
\|\big(J|_{\R\times S_{\rho,r}}\big)'(\s,u)\|:=\sup_{(\gamma,v)\in T_{(\s,u)}(\R\times S_{\rho,r}) \quad \|(\gamma,v)\|_{\R\times\X}=1} \langle J'(\s,u), (\gamma,v)\rangle_{(\R\times \X)^{-1}, (\R\times \X)},
\end{equation}
where the tangent space of $\R\times S_{\rho,r}$ at $(\s,u)$ is $T_{(\s,u)}(\R\times S_{\rho,r}):=\R\times T_u (S_{\rho,r})$ with
$$T_u (S_{\rho,r})=\{v\in \X_r : \langle u,v \rangle_{L^2}=0\}.$$

As a consequence of Ekeland's variational principle, following \cite{Mwlm}, we get the result below.

\begin{lemma}\label{le:ekeland}
Let $\eps>0$. Suppose that $(\tilde\s,\tilde\gamma)\in \Sigma \times\Gamma$ satisfies 
\begin{equation}\label{2.6w}
\max_{t \in [0,1]}J\big(\tilde\s(t),\tilde\gamma(t)\big)\le m+\eps,
\end{equation}
then there exists $(\s, u)\in \R\times S_{\rho,r}$ such that
\begin{enumerate}[label=(\arabic{*}),ref=\arabic{*}]
\item \label{eke1} ${\rm dist}_{\R \times S_{\rho,r}}\big((\s,u),(\tilde\s([0,1]),\tilde\gamma([0,1]))\big)\le 2 \sqrt{\eps}$;
\item \label{eke2} $J(\s,u)\in [m-2\eps,m+2\eps]$;
\item \label{eke3} $\|\big(J|_{\R\times S_{\rho,r}}\big)'(\s,u)\| \leq 8 \sqrt{\eps}$.
\end{enumerate}
\end{lemma}
\begin{proof}
By contradiction, suppose that there exists $(\tilde\s,\tilde\gamma)\in \Sigma \times\Gamma$ such that \eqref{2.6w} holds but \eqref{eke1}, \eqref{eke2} or \eqref{eke3} are not satisfied, namely for any $(\s, u)\in \R\times S_{\rho,r}$ at least one of the following conditions is true
\begin{enumerate}[label=($\arabic{*}')$,ref=$\arabic{*}'$]
\item \label{eke1'} ${\rm dist}_{\R \times S_{\rho,r}}\big((\s,u),(\tilde\s([0,1]),\tilde\gamma([0,1]))\big)> 2 \sqrt{\eps}$;
\item\label{eke2'}  $J(\s,u)\notin [m-2\eps,m+2\eps]$;
\item\label{eke3'} $\|\big(J|_{\R\times S_{\rho,r}}\big)'(\s,u)\| > 8 \sqrt{\eps}$.
\end{enumerate}
Then \eqref{eke3'} holds for any $(\sigma, u)$ such that \eqref{eke1} and \eqref{eke2} are satisfied.
Thus, applying Lemma 5.15 in \cite{Mwlm} with $X=\R\times\X_r$, $V=\R\times S_{\rho,r}$, $\varphi=J$,  $S=(\tilde\s([0,1]),\tilde\gamma([0,1]))\subset V$, $\delta=\sqrt{\varepsilon}$, there exists $\eta\in \C\big([0,1]\times(\R\times S_{\rho,r}), \R\times S_{\rho,r}\big)$ such that
\begin{equation}\label{5.151}
\eta\big(t, (\s, u)\big)=(\s, u), \quad\text{if }t=0\text{ or if } (\s,u)\notin J^{-1}([m-2\varepsilon,m+2\varepsilon])\cap \tilde S_{2\sqrt{\varepsilon}},
\end{equation}
\begin{equation}\label{5.152}
\eta\big(1,J^{m+\varepsilon}\cap  \tilde S_{2\sqrt{\varepsilon}}\big)\subset J^{m-\varepsilon},
\end{equation}
where $J^a:=\{v\in \R\times S_{\rho,r} : J(v)\le a\}$ and $ \tilde S_{\delta}:=\{v\in \R\times S_{\rho,r} : dist_{\R\times \X}(v, \tilde{S})\le \delta\}$.
Moreover
$J\big(\eta(\cdot, (\s, u))\big)$ is nonincreasing for every $(\s, u)\in \R\times S_{\rho,r},$
\\
By \eqref{muou1}, we can assume that
$$m-2\varepsilon>\max\{J(0,u_0),J(0,u_1)\}.$$
We define $\beta(t):=\eta\big(1,(\tilde\s(t),\tilde\gamma(t))\big)$, then by using \eqref{2.6w}, \eqref{5.151} and \eqref{5.152} we have 
$$\beta(0):=\eta\big(1,(\tilde\s(0),\tilde\gamma(0))\big)=\big(\tilde\s(0),\tilde\gamma(0)\big)=(0,u_0),$$
$$\beta(1):=\eta\big(1,(\tilde\s(1),\tilde\gamma(1))\big)=\big(\tilde\s(1),\tilde\gamma(1)\big) =(0,u_1),$$
so that $\beta\in \Sigma\times\Gamma$. Then, from \eqref{2.6w} and \eqref{5.152} we get
$$m\le \sup_{t\in[0,1]}J\big(\beta(t)\big)=\sup_{t\in[0,1]}J\left(\eta\big(1,(\tilde\s(t),\tilde\g(t))\big)\right)\le m-\eps$$
which is a contradiction.
\end{proof}

With all these results in our hands, it is easy to prove the existence of a (PS) sequence for $J|_{ \R \times S_{\rho,r}}$ at the mountain pass level $m$ verifying other crucial properties.

\begin{proposition}\label{pr:sequence}
There exists a sequence $\{(\s_n,u_n)\} \subset \R \times S_{\rho,r}$ such that, as $n \to +\infty$, we get 
	\begin{enumerate}[label=(\arabic{*}),ref=\arabic{*}]
		\item \label{seq1}$\s_n \to 0$;
		\item \label{seq2}$J(\s_n,u_n)\to m$; 
		\item \label{seq3}$\de_\s J(\s_n,u_n)\to 0$; 
		\item \label{seq4}$\de_u J|_{\R \times S_{\rho,r}}(\s_n,u_n)\to 0$ in $\X'_r$. 
	\end{enumerate}
	
\end{proposition}
\begin{proof}
	In view of Lemma \ref{le:ekeland} and \eqref{rema}, we conclude by letting $\eps\to0$.
\end{proof}

Now it only remains to prove Theorem \ref{mainb2} whose statement is given in the Introduction.

\begin{proof}[Proof of Theorem \ref{mainb2}]
By Proposition \ref{pr:sequence}, there exists a (PS) sequence $\{(\s_n,u_n)\} \subset \R \times S_{\rho,r}$ at level $m$ with $\s_n \to 0$ and satisfying the following
\begin{equation}\label{sistemab}
\begin{cases}
\dis\frac{e^{-N\s_n}}2 \irn B\big(e^{\s_n(N+2)}|\n u_n|^2\big)\, dx-\frac{e^{\s_n N(p/2-1)}}{p}\|u_n\|_p^p=m+ o(1),
\\[7mm]
\dis\frac{Ne^{-N\s_n}}{2}\!\!\irn\! \!\!B\big(e^{\s_n(N+2)}|\n u_n|^2\big)\, dx
-\frac{N+2}{2}e^{2\s_n}\!\!\irn \!\!b\big(e^{\s_n(N+2)}|\n u_n|^2\big)|\n u_n|^2\, dx
\\[2mm]
\dis \hspace{5.5cm}
+\frac{N(p-2)}{2p}e^{\s_n N(p/2-1)}\|u_n\|_p^p = o(1),
\\[7mm]
\dis e^{2\s_n}\!\irn \!b\big(e^{\s_n(N+2)}|\n u_n|^2\big)|\n u_n|^2\, dx-e^{\s_nN(p/2-1)}\|u_n\|_p^p +\rho^2 \l_n = o(1)\|u_n\|_{\X},
\end{cases}
\end{equation}
where $\l_n$ is a Lagrange multiplier. 
From the first and the second equations of the system, by \eqref{growthb1} and \eqref{b3} we get
\begin{align*}
C+ o(1)&=\frac{Npe^{-N\s_n}}{4}\!\!\irn \!\! B\big(e^{\s_n(N+2)}|\n u_n|^2\big) dx
-\frac{N+2}{2}e^{2\s_n}\!\!\irn \!\! b\big(e^{\s_n(N+2)}|\n u_n|^2\big)|\n u_n|^2 dx\\
&\ge
\left(\frac{Npe^{-N\s_n}}{4}\a -\frac{N+2}{2}e^{2\s_n}\right)\irn b\big(e^{\s_n(N+2)}|\n u_n|^2\big)|\n u_n|^2\, dx\\
&\ge C'(\|\n u_n\|_2^2+\|\n u_n\|_q^q),
\end{align*}
where $C,C'>0$ are constants.
Therefore, since $\s_n \to 0$, as $n \to+\infty$, by \eqref{equ2.2} we deduce that $\{u_n\}$ is a bounded sequence in $D^{1,2}(\RN)\cap D^{1,q}(\RN)$ and so also  in $\X_r$, being $u_n\in S_{\rho,r}$. By the third equation of \eqref{sistemab} and \eqref{growthb1}, also $\{\l_n\}$ is bounded in $\R$.
So there exist $u \in \X_r$ and $\l\in \R$ such that
$$ u_n \rightharpoonup u\text{ in } \X_r,\quad u_n\to u\hbox{ in }L^p(\RN), \quad
u_n \to u\,\,\text{ a.e in }\,{\R^N}, \quad \l_n\to \l,$$
up to subsequences, as $n \to +\infty$.
Since $I'(u_n)+\l_n u_n \to 0$ in $\X_r'$, then arguing as in the proof of Lemma \ref{LeMa3.4.1}, 
we deduce that 
\[
\n u_n \to \n u \quad \text{in }L^2_{loc}(\RN)\cap L^q_{loc}(\RN), \qquad \n u_n \to \n u \quad \text{a.e. in }\RN.
\]
Therefore, since $I'(u)+\l u = 0$ in $\X_r'$, we deduce that 
\begin{equation}\label{cacca}
\irn b(|\n u|^2)|\n u|^2\, dx+\l\|u\|_2^2=\|u\|_p^p.
\end{equation}
{\em Step 1:} $\l\neq 0$.
\\
Assume by contradiction that $\l=0$, then, by \eqref{cacca},  
\[	\irn b(|\n u|^2)|\n u|^2\, dx=\|u\|_p^p.\]
By the weak lower semicontinuity and by the third equation of \eqref{sistemab} and the strong convergence in $L^p(\RN)$,
\begin{align*}
\irn b(|\n u|^2)|\n u|^2\, dx&\le
\liminf_{n\to +\infty}e^{2\s_n}\irn b(e^{\s_n(N+2)}|\n u_n|^2)|\n u_n|^2\, dx\\
&= \dis\liminf_{n\to +\infty} \left[e^{N\s_n(p/2-1)}\|u_n\|_p^p + o(1)\right]
\\
&=\|u\|_p^p
=\irn b(|\n u|^2)|\n u|^2\, dx
\end{align*}
and so, up to a subsequence,
\begin{align}
\irn b(|\n u|^2)|\n u|^2\, dx 
&= \lim_{n \to +\infty}\irn b(e^{\s_n(N+2)}|\n u_n|^2)|\n u_n|^2\, dx. \label{conv1b}
\end{align}
In view of \cite[Lemma 2.5]{MP} equation \eqref{conv1b} implies that $u_n \to u$ strongly in $D^{1,2}(\R^N)\cap D^{1,q}(\R^N)$. So, by the first equation of \eqref{sistemab}, we have $I(u)=m\neq0$, and so also $u\neq 0$.
\\
Since we are assuming, by contradiction, that $\lambda_n\to 0$, then
$$\irn b\big(e^{\s_n(N+2)}|\n u_n|^2\big)|\n u_n|^2\, dx-\|u_n\|_p^p = o(1),$$
which, together with the second equation in \eqref{sistemab}, gives
$$\Big(\frac{N+2}{2}-\frac{N(p-2)}{2p}\Big)\irn b \big(e^{\s_n(N+2)}|\n u_n|^2\big)|\n u_n|^2\,dx=\frac{N}{2}\irn B \big(e^{\s_n(N+2)}|\n u_n|^2\big)\, dx+ o(1).$$
Now, by using that $u_n \to u$ strongly in $D^{1,2}(\R^N)\cap D^{1,q}(\R^N)$, we have
$$\Big(\frac{N+2}{N}-\frac{N(p-2)}{Np}\Big)\irn b \big(|\n u|^2\big)|\n u|^2\,dx=\irn B \big(|\n u|^2\big)\, dx,$$
that is
$$\irn\left[\frac{2(N+p)}{Np} b \big(|\n u|^2\big)|\n u|^2- B \big(|\n u|^2\big)\right]\, dx=0,$$
which is in contradiction with \eqref{b3}, since $u\neq 0$.
\\
{\em Step 2:} $\l> 0$.
\\
Combining the second and the third equation of \eqref{sistemab} and using \eqref{b3}, we have
\begin{align*}
N(p-2)\rho^2 \l_n +& o(1)
=-N(p-2)\irn b\big(e^{\s_n(N+2)}|\n u_n|^2\big)|\n u_n|^2\, dx
+N(p-2)\|u_n\|_p^p 
\\
&=\irn \left[2(N+p)b\big(e^{\s_n(N+2)}|\n u_n|^2\big)|\n u_n|^2-Np B\big(e^{\s_n(N+2)}|\n u_n|^2\big) \right]\, dx\ge 0.
\end{align*}
So, since $\l_n \to \l\neq 0$, we conclude that  $\l> 0$.
\\
{\em Step 3:} $u_n \to u$ strongly in $\X_r$.
\\
Since  $u_n \weakto u$ in $\X_r$, by the weak lower semicontinuity we have that
\begin{align}
\irn b(|\n u|^2)|\n u|^2\, dx 
&\le \liminf_{n \to +\infty}\irn b(e^{\s_n(N+2)}|\n u_n|^2)|\n u_n|^2\, dx,\label{caccab}
\\
\|u\|_2^2 &\le \liminf_{n \to +\infty}\|u_n\|_2^2=\rho^2.\label{caccac}
\end{align}
Therefore, by the third equation in \eqref{sistemab} and \eqref{cacca}, we have
\begin{align*}
\irn b(|\n u|^2)|\n u|^2\, dx&+\l\|u\|_2^2\le
\liminf_{n\to +\infty}\left[e^{2\s_n}\irn b(e^{\s_n(N+2)}|\n u_n|^2)|\n u_n|^2\, dx
+\l_n \rho^2\right]\\
&= \dis\liminf_{n\to +\infty} \left[e^{N\s_n(p/2-1)}\|u_n\|_p^p + o(1)\right]
\\
&=\|u\|_p^p
=\irn b(|\n u|^2)|\n u|^2\, dx+\l \|u\|_2^2.
\end{align*}
So, since $\l>0$, by \eqref{caccab} and \eqref{caccac}, we deduce that, up to a subsequence,
\begin{align}
\irn b(|\n u|^2)|\n u|^2\, dx 
&= \lim_{n \to +\infty}\irn b(e^{\s_n(N+2)}|\n u_n|^2)|\n u_n|^2\, dx, \label{conv3b}
\\
\|u\|_2^2 &= \lim_{n \to +\infty}\|u_n\|_2^2=\rho^2.\label{conv3c}
\end{align}
In view of \cite[Lemma 2.5]{MP} equation \eqref{conv3b} implies that $u_n \to u$ strongly in $D^{1,2}(\R^N)\cap D^{1,q}(\R^N)$, while \eqref{conv3c} gives the strong convergence also in $L^2(\RN)$.
This assures that $u\in S_{\rho,r}$ is a radial normalized solution to \eqref{Pb} concluding the proof.

\end{proof}

\begin{remark}
Observe that the following inequalities 
$$\big(1+\frac{2}{N}\big)q<p<2^*\hbox{ and }q>N$$ hold for $N=3$. Therefore, there is the attempt to consider also the Born-Infeld problem in the $L^2$-supercritical case, in the three-dimensional case. The idea would be to apply Theorem \ref{mainb}  to the penalization $a_\t$ and then to get rid of the penalization with a suitable uniform estimate on the $L^\infty$-norm of the solution, as done in Section \ref{res}. However  \eqref{b3} is not satisfied by $a_\t$, for $\t$ near zero, and this assumption seems to be crucial in order to prove the boundedness of the {\em special} (PS) sequence obtained by Proposition \ref{pr:sequence}.
So the $L^2$-supercritical case for the Born-Infeld problem still remains an open problem.
\end{remark}

\section{Problem \eqref{Pb} in $L^2$- and $L^q$-critical cases}\label{crit}
%
%
%
%
%

We start considering the $L^2$-critical case. We assume, therefore, $p=2+4/N$, and fix $q>1$. Moreover we define $0<\rho_*<\rho^*$ as follows 
$$\rho_*:=c_1^{N/4}\|W_{p}\|_{2}\qquad \rho^*:=c_2^{N/4}\|W_{p}\|_{2}$$ 
where $c_1$ and $c_2$ are introduced in \eqref{b1} and $W_{p}$ is  defined in Lemma \ref{lem2.1}. 

\begin{proof}[Proof of Theorem \ref{th1.2}]

\eqref{i-th1.2} We first observe that $p\delta_{p}=2$ if $p=2+4/N$. For any fixed $u \in S_{\rho}$, the Gagliardo-Nirenberg inequality \eqref{equ2.2} and \eqref{energyFinb} indicate that
\begin{align} \label{4-lower-bounded1.1}
I(u)\ge\frac{c_1}{2} \|\nabla u\|_{2}^2+\frac{c_1}{2} \|\nabla u\|_{q}^q-\frac{1}{p} {\|u\|}_{ p }^{ p }
\geq \biggl(\frac{c_1}{2}-\frac{\mathcal{C}^{ p }_{N,p}}{p} \rho^{p-2}\biggr) \left\|\nabla u\right\|_{2}^{2}+\frac{c_1}{2} \|\nabla u\|_{q}^q>0
\end{align} 
when $\rho \in (0,c_1^{N/4} \|W_{p}\|_{2}]$. So we have $m(\rho)\geq0$. In addition, we have $u_{t}(x)=t^{\frac{N}{2}}u(tx)\in S_{\rho}$ and hence
\begin{equation*}\label{upperbound}
m(\rho)\leq I(u_t)\leq\frac{c_2 t^{2}}{2} \|\nabla u\|_{2}^2+\frac{c_2 t^{q(1+\delta_{q})} }{2} \|\nabla u\|_{q}^q-\frac{t^{ 2}}{ p } {\|u\|}_{ p }^{ p } \to 0
\end{equation*}
as $t\to0$. This implies that $m(\rho)=0$. From \eqref{4-lower-bounded1.1}, we also deduce that $m(\rho)$ can not be attained.

\eqref{ii-th1.2} First, since $W_p$ is a solution of \eqref{eqlem2.1}, we get
\begin{equation*} \label{TeStfunC3.3}
 \| \nabla W_p \|_{2}^{2}=\| W_p \|_{2}^{2}=\frac{2}{p} \| W_p \|_{p}^{p}.
\end{equation*}
In addition, it follows from \cite[Proposition 4.1]{Gdas} that
\begin{equation*} \label{TeStfunC3.4}
W_p(x),|\nabla W_p(x)|=O\left(|x|^{-\frac{N-1}{2}} e^{-|x|}\right)  \quad \text { as }|x| \rightarrow +\infty.
\end{equation*}
 Let $\varphi(x) \in \cC_{0}^{\infty}\left(\mathbb{R}^{N}\right)$ be a nonnegative radial function such that $\varphi(x)=1$ if $|x| \leq 1$ and $\varphi(x)=0$ if $|x| \geq 2$. For $\tau>0$ and $R>0$, we denote
$$\phi_{1}(x)=A_{\tau, R} \frac{ (\tau \rho)^{ \frac{N}{2} }  }{\|W_p\|_{2}} \varphi \Big( \frac{x}{R} \Big) W_p ( \tau x ) ,$$
where $A_{\tau, R}>0$ is chosen such that $\int_{\mathbb{R}^{N}} \phi_{1}^{2}\, dx=\rho^2$.
As in \cite{BY} we have
\begin{align*}
\|\nabla \phi_{1}\|^{2}_2
&= \tau^{2} \frac{ \rho^{2} }{\|W_p\|_{2}^{2}}  \|\nabla W_p\|^{2}_2  +o(1),
\\
\|\nabla \phi_{1}\|^{q}_q
 &\leq  \tau^{ q(1+\delta_{q}) } \frac{  \rho^q}{\|W_p\|_{2}^{q}}  \|\nabla W_p\|^{q}_q  +o(1),
\\
\|\phi_{1}\|^{p}_p
 &= \tau^{ p\delta_{p} } \frac{\rho^p}{\|W_p\|_{2}^{p}}\| W_p\|^{p}_p +o(1),
\end{align*}
as $ \tau R \rightarrow +\infty$

Therefore, letting $\tau \to +\infty$, 
\begin{equation}\label{lowerbound-4.9}\begin{aligned}
m(\rho)& \leq  I(\phi_{1}) \leq \frac{c_2\tau^{2}}{2} \frac{ \rho^2}{\|W_p\|_{2}^{2}}  \|\nabla W_p\|^{2}_2
+\frac{c_2 \tau^{q(1+\delta_{q})} }{2} \frac{ \rho^q}{\|W_p\|_{2}^{q}}  \|\nabla W_p\|^{q}_q  \\
&\qquad-\frac{\tau^{ p\delta_{p} }}{p}  \frac{\rho^p}{\|W_p\|_{2}^{p}}\| W_p\|^{p}_p  +o(1), \\
&= \tau^{2}\rho^2  \biggl(  \frac{c_2 \|\nabla W_p\|_2^{2}}{2\|W_p\|_{2}^{2}}  -   \frac{ \rho^{p-2} \| W_p\|_p^{p}  }{p\|W_p\|_{2}^{p}}  \biggr)  +\frac{c_2\tau^{ q(1+\delta_{q}) }}{2} \frac{  \rho^q}{\|W_p\|_{2}^{q}}  \|\nabla W_p\|_{q}^{q} +o(1).
\end{aligned}\end{equation}
Since $1<q<2$ and $\rho>c_2^{N/4}\|W_{p}\|_{2}$, we have
$$ \frac{c_2 \|\nabla W_{p}\|_2^{2}}{2\|W_{p}\|_{2}^{2}}  -   \frac{ \rho^{{p}-2} \|W_{p}\|_{p}^{{p}}  }{{p}\|W_{p}\|_{2}^{{p}}}<0\text{ and } q(1+\delta_{q})<2 .$$
Letting $\tau \to +\infty$ in \eqref{lowerbound-4.9}, we have $m(\rho) =-\infty$.
\end{proof}

\begin{remark}
Since $\rho_*<\rho^*$ for $c_1<c_2$, when $1<q<2$, we cannot prove nonexistence for any $\rho>0$, as in \cite{BY}.
\end{remark}

Now we conclude considering the $L^q$-critical case. Namely we consider $p=q+2q/N$.
Let us define
$$
\hat \rho_{*} :=\biggl(\frac{c_1q(N+2)}{2N\mathcal{K}_{N,p}^{q(N+2)/N}}\biggr)^{N/2q},
 \qquad 
\hat{\rho}^{*}
:=\biggl[\frac{ c_2\|\nabla W_p\|_q^q}{\|W_p\|_2^{2(N-q)/N}}\biggr]^{N/2q},$$
where $\mathcal{K}_{N,p}$ has been introduced in Lemma \ref{lem2.2}.

\begin{proof}[Proof of Theorem \ref{th1.2bis}]
	\eqref{i-th1.2bis} Following proof of Theorem \ref{th1.2} by using \eqref{pq-equ2.2} in place of \eqref{equ2.2} and observing that $p\nu_{p,q}=q$ if $p\!=\!q+2q/N$, for $u\in S_\rho$ we have     
	\begin{align*}
	I(u)&\ge\frac{c_1}{2} \|\nabla u\|_{2}^2+\frac{c_1}{2} \|\nabla u\|_{q}^q-\frac{1}{ {p}} {\|u\|}_{  {p} }^{  {p} }  \\
	&\geq \frac{c_1}{2} \|\nabla u\|_{2}^2+\frac{c_1}{2} \|\nabla u\|_{q}^q-\frac{\mathcal{K}^p_{N,p}}{ {p}} \left\|\nabla u\right\|_{q}^{p\nu_{p,q} } \rho^{p(1-\nu_{p,q} )}  \\
	&=\frac{c_1}{2} \|\nabla u\|_{2}^2+ \Big(\frac{c_1}{2}-\frac{\mathcal{K}^p_{N,p}}{ {p}}  \rho^{p-q}  \Big) \left\|\nabla u\right\|_{q}^{q} >0
	\end{align*}
	when $\frac{\mathcal{K}^p_{N,p}}{ {p}}  \rho^{p-q}\le\frac{c_1}{2}$, that is if $\rho\le \hat \rho_{*}$. So that, taking into account \eqref{upperbound}, if $0<\rho\le \hat \rho_{*}$, then $m(\rho)=0$ and there is no minimizer of $m(\rho)$. Thus, we conclude.

	\eqref{ii-th1.2bis} In a similar way, since $p\delta_p=q(1+\delta_q)$ then \eqref{lowerbound-4.9} becomes
	$$m(\rho)\leq \frac{c_2\tau^{2}}{2} \frac{ \rho^2}{\|W_p\|_{2}^{2}}  \|\nabla W_p\|^{2}_2+\tau^{ p\delta_{p} }\rho^q\biggl[\frac{c_2 \|\nabla W_p\|_q^q}{2\|W_p\|_{2}^{q}}-\frac{\rho^{p-q}\|W_p\|_p^p}{p\|W_p\|_{2}^{p}}\biggr]+o(1),$$  
	as $\tau R\to+\infty$. Note that $2\!<\!q\!<\!N$, so we have $ p\delta_{p}=q(1+\delta_{q})>2$. In view of this fact and $\rho>\hat{\rho}^{*}$, we let $\tau \to +\infty$ and obtain $m(\rho) =-\infty$.  
\end{proof}

\begin{remark}
Note that, we get the nonexistence for $\rho$ small for every $q>1$ in the $L^2$-critical case, while in the $L^q$-critical case only for $1<q<N$.
\end{remark}

\section*{Acknowledgments}
L. Baldelli is member of the {\em Gruppo Nazionale per l'Analisi Ma\-te\-ma\-ti\-ca, la Probabilit\`a e le loro Applicazioni}
(GNAMPA) of the {\em Istituto Nazionale di Alta Matematica} (INdAM).
L. Baldelli was partly supported by National Science Centre, Poland (Grant No. 2020/37/B/ST1/02742), by INdAM-GNAMPA Project 2023 titled {\em Problemi ellittici e parabolici con termini di reazione singolari e convettivi} (E53C22001930001) and by the IMAG-Maria de Maeztu Excellence Grant CEX2020-001105-M funded by MICINN/AEI. \newline J. Mederski was partly supported by National Science Centre, Poland (Grant No. 2017/26/\\E/ST1/00817).\\
A. Pomponio is partly supported by PRIN projects 2017JPCAPN {\em Qualitative and quantitative aspects of nonlinear}
PDEs, P2022YFAJH {\em Linear and Nonlinear PDEs: New directions and Applications}, and by INdAM - GNAMPA Project 2022  {\em Modelli EDP nello studio problemi della fisica moderna} and INdAM - GNAMPA Project 2023 
{\em Metodi variazionali per alcune equazioni
di tipo Choquard}.

\end{document}